\documentclass[10 pt,a4paper,leqno]{article}
\usepackage{amsmath}

\RequirePackage{amsmath}
\RequirePackage{amssymb}
\RequirePackage{ifthen}
\RequirePackage{theorem}
\RequirePackage{pslatex}

\title{\sc The $q$-analogue of the wild fundamental group (I)}

\author{J.-P. Ramis\footnote{Institut de France (Acad\'emie des Sciences)},\;
J. Sauloy
\footnote{Laboratoire Emile Picard, CNRS UMR 5580, U.F.R. M.I.G.,
118, route de Narbonne, 31062 Toulouse CEDEX 4}}

\date{}

\def\sq{\sigma_q}

\def\C{{\mathbf C}}
\def\Q{{\mathbf Q}}
\def\Z{{\mathbf Z}}
\def\R{{\mathbf R}}
\def\N{{\mathbf N}}
\def\U{{\mathbf U}}

\def\F{{\mathcal{F}}}

\def\Bd{{\mathcal{B}_{q}^{(\delta)}}}

\def\O{{\mathcal{O}}}
\def\M{{\mathcal{M}}}
\def\Ma{{\text{Mat}}}

\def\D{{\mathcal{D}_{q,K}}}
\def\gr{{\text{gr}}}

\def\Eq{{\mathbf{E}_{q}}}
\def\EE{{\mathcal{E}}}
\def\G{{\mathfrak{G}}}
\def\g{{\mathfrak{g}}}

\def\Ra{{\C\{z\}}}
\def\Ka{{\C(\{z\})}}
\def\Raqd{{\C\{\xi\}_{q,\delta}}}
\def\Kaqd{{\C(\{\xi\})_{q,\delta}}}
\def\Raq(d){{\C\{\xi\}_{q,(\delta)}}}
\def\Kaq(d){{\C(\{\xi\})_{q,(\delta)}}}
\def\Ree{{\mathcal{O}(\C)}}
\def\Ke{{\mathcal{O}(\C)[z^{-1}]}}
\def\Rf{{\C[[z]]}}
\def\Kf{{\C((z))}}
\def\Rw{{\mathcal{O}(\C^{*})}}
\def\Kw{{\mathcal{M}(\C^{*})}}
\def\Rwg{{\mathcal{O}(\C^{*},0)}}
\def\Kwg{{\mathcal{M}(\C^{*},0)}}
\def\DM{{DiffMod(K,\sq)}}

\def\Der{{\dot{\Delta}}}
\def\Derc{{\Der_{\overline{c}}}}
\def\Derdc{{\Der^{(\delta)}_{\overline{c}}}}
\def\Derdx{{\Der^{(\delta)}_{\overline{\xi}}}}
\def\Derdxab{{\Der^{(\delta,\alpha,\beta)}_{\overline{\xi}}}}
\def\Lsca{{LS_{\overline{c},a}}}

\def\l{\left}
\def\r{\right}

\def\th{\theta}

\def\Im{\text{Im}}
\def\Ker{\text{Ker}}
\def\Id{\text{Id}}

\def\limproj{\mathop{\oalign{lim\cr \hidewidth$\longleftarrow$\hidewidth\cr}}}

\def\St{\mathfrak{St}}
\def\st{\mathfrak{st}}

\newtheorem{thm}{Theorem}[section]
\newtheorem{lemma}[thm]{Lemma}
\newtheorem{prop}[thm]{Proposition}
\newtheorem{cor}[thm]{Corollary}

\def \Pr {\textsl{Proof. - }}

\begin{document}

\maketitle

\bigskip \hrule \bigskip

\begin{abstract}
{\small 
We describe an explicit construction of galoisian Stokes operators
for irregular linear q-difference equations. These operators are
parameterized by the points of an elliptic curve minus a finite
set of singularities. Taking residues at these singularities, one
gets q-analogues of alien derivations which freely generate the
Lie algebra of the Stokes subgroup of the Galois group.}
\end{abstract}

\bigskip \hrule 



\section{Introduction}

In this paper we return to the local analytic classification 
of $q$-difference modules. In \cite{RSZ1} we gave such a classification 
in Birkhoff style \cite{Birkhoff1, Birkhoff2}, using normal forms and index theorems. In \cite{JSSTO} 
(cf. also \cite{RSZ}, \cite{RSZ2}) appears a classification using 
non abelian cohomology 
of sheaves on an elliptic curve. Here our aim is to give a classification 
based upon a ``fundamental group'' and its finite dimensional representations, 
in the style of the Riemann-Hilbert correspondence. At the abstract level 
there exists such a classification: the fundamental group is the tannakian 
Galois group of the tannakian category of our $q$-modules, but we want more 
information: our final aim is to get a \emph{smaller} fundamental group 
(as small as possible !) which is Zariski dense in the tannakian Galois 
group and to describe it \emph{explicitly}. (As a byproduct, we shall get 
a complete description of the tannakian Galois group itself.) It is 
important to notice that the tannakian Galois group is an algebraic 
object but that the construction of the smaller group is based upon 
\emph{transcendental techniques} (topology and complex analysis). \\
This program is only partially achieved in this paper: even if we have 
the main ingredients, part of the work remains conjectural. As in the 
differential case, the construction of the fundamental group is a 
Russian-dolls construction using semi-direct products (heuristically, 
going from interior to exterior,  each new ``slice of infinitesimal 
neighborhood'' of the origin (each ``scale'') corresponds to an invariant 
subgroup in a new semi-direct product). At the end there is a fascinating 
parallel between the differential and the $q$-difference case. However, 
it has been impossible (for us) to mimick the differential approach 
(essentially based upon the concept of \emph{solutions}); instead, 
we shall follow a new path. Before the description of our approach 
and of our results in the $q$-difference case, let us recall, without 
proof, what happens in the differential case (in this case the ``fundamental 
group'' is the {\it wild fundamental group} introduced by the first author). 
We shall not use these results but we need them at the heuristic level 
in order to develop our $q$-analogies, therefore we shall insist on 
the underlying geometric ideas. \\

We shall use \emph{tannakian categories} as an essential tool (cf. \cite{DF}, 
\cite{Simpson}, 6, page 67). First recall some basic facts. Assume 
$\mathcal E$ is a \emph{neutral tannakian category}, with fiber functor 
$\omega$ (to the category of $\C$-vector spaces), then 
$Aut^{\otimes}(\omega )$ has a structure of \emph{complex pro-algebraic 
affine group scheme}; we shall call it the tannakian group of the tannakian 
category $\mathcal E$. The category $\mathcal E$ is isomorphic to 
the category of finite dimensional representations of $Aut^{\otimes}(\omega )$ 
(by \emph{definition}, such a representation factors through a representation 
of one of the algebraic quotients). Conversely, if $G$ is a complex 
pro-algebraic group, its category of complex representations $Rep_{\C}(G)$ 
is a neutral tannakian category with a natural fiber functor (the obvious
forgetful functor) $\omega_G$ and $G=Aut^{\otimes}(\omega_G,Rep_{\C}(G))$; 
the complex space ${\bf aut}^{\otimes}(\omega_G,Rep_{\C}(G))$ of Lie-like 
$\otimes$-endomorphisms of the fiber functor $\omega_G$ is the Lie-algebra 
of $Aut^{\otimes}(\omega_G,Rep_{\C}(G))$. 

Assuming that $\Gamma$ is a \emph{finitely generated} group,
a \emph{pro-algebraic completion} of $\Gamma$ is, by definition, 
a \emph{universal} pair $(\iota_{al},\Gamma^{al})$ where 
$\iota_{al}: \Gamma \rightarrow \Gamma^{al}$ is a group homomorphism 
from $\Gamma$ to a {\it pro-algebraic group} $\Gamma^{al}$.  It is unique 
up to an isomorphism of proalgebraic groups. A finite dimensional 
representation of $\Gamma$ clearly factors through a proalgebraic 
completion of $\Gamma$. We can get a pro-algebraic completion of 
$\Gamma$ from the tannakian mechanism: $Rep_{\C}(\Gamma )$ is a neutral 
tannakian category with a natural fiber functor (the obvious forgetful 
functor) $\omega$, and the group $G=Aut^{\otimes}(\omega ,Rep_{\C}(\Gamma ))$ 
is a pro-algebraic completion of $\Gamma$. The groups $\Gamma$ and $G$ 
have the same representations, more precisely the natural homomorphism 
of groups $\Gamma \rightarrow G$ induces an isomorphism of tannakian 
categories: $Rep_{\C}(G) \rightarrow Rep_{\C}(\Gamma )$. We shall encounter 
below similar situations associated to different classification problems 
(in a little more general setting: $\Gamma$ will not be in general 
a finitely generated group). \\
The first example (our baby example) is the category of local meromorphic 
regular-singular connections, or equivalently the category 
${\mathcal D}^{(0)}_f$ of regular singular $\mathcal D$-modules, 
where $\mathcal D=\Ka[d/dz]$ ($\Ka$ is the field of fractions of $\Ra$). 
A meromorphic connection is equivalent to an equivalence class
of differential systems $\Delta_A:~\frac{dY}{dx}=AY$ up to the 
gauge-equivalence: $\Delta_A\sim \Delta_B$ if and only if there exists 
$P\in Gl_n(\Ka)$ such that $B=P^{-1}AP-P^{-1}dP/dx$. We consider the 
fundamental group $\pi_1(D^*,d)$ of a germ at zero  of punctured disk, 
pointed on a germ of direction $d$. We choose a generator $\gamma$ 
(a one turn loop in the positive sense) and we get an isomorphism 
$\Z\rightarrow \pi_1(D^*,d)$, $n\mapsto \gamma^n$.  Then, by a very 
simple application of the Riemann-Hilbert correspondance, our category 
${\mathcal D}^{(0)}_f$ is equivalent (via the \emph{monodromy 
representation}) to the category of finite dimensional representations 
of the fundamental group $\pi_1(D^*,d)$. A ${\mathcal D}$-module $M$ 
corresponds to a representation $\rho_M$ of $\pi_1(D^*,d)$. \\
We can apply the tannakian machinery to the group $\Gamma =\Z$ 
(or equivalently to $\Gamma=\pi_1(D^*,d)$). Then our category 
${\mathcal D}^{(0)}_f$ is equivalent to the category of representations 
of $\pi_1(D^*,d)$: a $\mathcal D$-module $M$  ``is'' a representation 
$\rho_M$ of the topological fundamental group $\pi_1(D^*,d)$, it is also 
equivalent to the category of representations of the pro-algebraic 
completion  $\pi_1^{\otimes}(D^*,d)$ of $\pi_1(D^*,d)$: 
a $\mathcal D$-module ``is'' a representation $\rho_M^{\otimes}$ of 
the tannakian fundamental group $\pi_1^{\otimes}(D^*,d)$. The ``small 
fundamental group'' is the topological group $\pi_1(D^*,d)$, the ``big 
fundamental group'' is the pro-algebraic group $\pi_1^{\otimes}(D^*,d)$. 
The small group is \emph{Zariski-dense} in the big group: the image 
of $\rho_M$ is the monodromy group of $M$, it is Zariski-dense in the image 
of $\rho_M^{\otimes}$ which ``is'' the differential galois group of $M$. \\
It is not difficult but important to understand the classification mechanism 
on this baby example: all the information is hidden in the group $\Z$ and 
we must extract it. The essential point is to understand the structure of 
the pro-algebraic completion of $\Z$. We can use the tannakian machinery
(this is ``folklore knowledge'', a reference is \cite{JSGAL}): 
the pro-algebraic hull of $\Z$ is $\Z^{al}={Aut}^{\otimes}(\omega)$, 
it is commutative and the product of its semi-simple part $\Z^{al}_s$ 
and its unipotent part  $\Z^{al}_u$:
$\Z^{al}_s=Hom_{gr}(\C^*,\C^*)$,  $\Z^{al}_u=\C$ (the additive group) and 
$\iota_{al}:\Z\rightarrow \Z^{al}$ is defined by $1\mapsto (id_{\C^*},1)$ 
($n\mapsto ((z\mapsto z^n),n)$). In order to understand what will happen 
in more difficult situations, it is interesting to understand the 
pro-algebraic completion of $\Z$ using regular singular differential 
equations and differential Galois theory. We shall start from  
$\pi_1(D^*,d)$ and shall ``compute'' its pro-algebraic completion 
using Riemann-Hilbert correspondance. The main tool is a universal 
Picard-Vessiot algebra ${\mathcal U}_f$  (\cite{Singer-vdPut2}). 
We consider some holomorphic functions on the Riemann surface of 
the logarithm: $\log~x$ and $x^{\alpha}=e^{\alpha \log~x}$ ($\alpha \in \C$). 
They generate over $\Ka$ a \emph{differential algebra}: 
 ${\mathcal U}_f=\Ka\big\{ (x^{\alpha})_{\alpha\in \C},\log~x\big\} $ 
(it is a \emph{simple} differential algebra; the brackets 
$\big\{\cdots\big\}$ mean ``differential algebra generated by ...). 
For each object $M$ 
of ${\mathcal D}^{(0)}_f$, the algebra ${\mathcal U}_f$ contains one 
and only one Picard-Vessiot algebra for $M$. Equivalently we can solve 
any regular-singular system 
$\Delta_A:~\frac{dY}{dx}=AY$ using ${\mathcal U}_f$ (that is we can find 
a fundamental matrix solution with entries in ${\mathcal U}_f$). 
The differential Galois group $G_f$ of ${\mathcal U}_f$ (or equivalently 
of its field of fractions) is a pro-algebraic group (${\mathcal U}_f$ is 
an inductive limit of finite type differential extensions). The monodromy, 
that is the action of the loop $\gamma$ is Galois, therefore we can 
identify $\gamma$ with an element of 
$G_f$ ($\gamma (x^{\alpha})=e^{2i\pi\alpha}x^{\alpha}$ 
and $\gamma (\log~x)=\log~x+2i\pi $), and we get an injective homomorphism 
of groups: $\iota :\pi_1(D^*,d) \rightarrow G_f$.
It is not difficult to prove that $(G_f,\iota )$ is a pro-algebraic 
completion of $\pi_1(D^*,d)$ and to compute $G_f$ (we shall admit 
these results here):  its semi-simple part 
$G_{f,s}=\check {\C/\Z}$ is the topological dual group of $\C/\Z$ 
considered as the inductive limit of its finitely generated subgroups;
its unipotent part $G_{f,u}$ is the differential Galois group of 
the extension $\Ka \big\{Ê\log\ z\big\} $, that is the additive 
group $\C$. We have an exact sequence of groups:
$$
0\rightarrow \Q/\Z \rightarrow \C /\Z \rightarrow \C/\Q \rightarrow 0
$$
and an exact sequence of dual groups:
$$
1\rightarrow {\bf T}_f\rightarrow G_{f,s} \rightarrow  
\hat {\Z}(1)\rightarrow 1.
$$
(The proalgebraic group ${\bf T}_f$ is the topological dual group of 
the group of ``monodromy exponents''). Here, 
$G_{f,s} =Hom_{gr}(\C/\Z,\C^*)\approx Hom_{gr}(\C^*,\C^*)$.  
The respective images of $\gamma$ in $G_{f,s}$ and in $G_{f,u}$ are
$z\mapsto z$ and $2i\pi$. \\

The next step is the study of the category of \emph{formal} connections, 
or equivalently the category ${\mathcal D}_{form}$ of 
$\hat{\mathcal D}$-modules, where $\hat{\mathcal D}=\Kf[d/dz]$ ($\Kf$ is the
field of fractions of  $\Rf$). It is a neutral tannakian category 
(cf. for instance \cite{Singer-vdPut2}). 
One associates to a $\hat{\mathcal D}$-module $M$ its 
\emph{Newton polygon} $N(M)$. The slopes of $N(M)$ are positive rational 
numbers. For sake of simplicity we shall limit ourselves to the full 
subcategory ${\mathcal D}_{form,int}$ of modules $M$ whose Newton polygon 
has only \emph{integer} slopes. In order to compute the corresponding 
``fundamental groups'', it is necessary to understand the formal 
classification of differential equations of order one: two such 
equations $dy/dx-\hat ay=0$ and
$dy/dx-\hat by=0$ ($\hat a,\hat b\in \Kf$) are formally equivalent 
if and only if $(b-a)dx\in d_{log}~\Kf$ 
($d_{log}~\Kf =\{ d\hat c/\hat c\vert~\hat c\in \Kf\} $). We have 
$d_{log}~\Kf =\Z \frac{dx}{x}\oplus \Rf dx$ 
and $\Kf dx/d_{log}~\Kf =\C/\Z \frac{dx}{x}\oplus \Kf/z^{-1}\Rf dz$. 
By integration $a \, dx\in \Kf dx$ gives
$q=\int a \, dx$ and we get an isomorphism between $\Kf dx/d_{log}~\Kf$ and
$\C/\Z~ \log~x\oplus \Kf/\Rf=\C/\Z~ \log~x\oplus \frac{1}{z}\C[\frac{1}{z}]$. 
To $\alpha \log~x \in\C/\Z~ \log~x$ corresponds 
$e^{\alpha \log~x}=x^{\alpha}$ (a solution of $dy/dx-\alpha y=0$), 
to $q\in \frac{1}{z}\C[\frac{1}{z}]$ corresponds $e^q$ (a solution 
of $dy/dx -q'y=0$). Therefore it is natural to introduce the differential 
algebra 
${\mathcal U}_{form,int}=
\Kf\big\{ (x^{\alpha})_{\alpha\in \C},
(e^q)_{q\in \frac{1}{z}\C[\frac{1}{z}]},\log~x\big\} $. 
It is a universal Picard-Vessiot algebra for the formal connections 
whose Newton polygons have only integer slopes and the differential 
Galois group of
${\mathcal U}_{form,int}$ is isomorphic to the tannakian Galois group 
of the category ${\mathcal D}_{form,int}$.
We consider $\frac{1}{z}\C[\frac{1}{z}]$ as a $\Z$-module. It has no 
torsion, it is an infinite dimensional lattice and we consider it as 
the inductive 
limit of its finite dimensional sublattices. The topological dual group 
of such a sublattice is a \emph{torus} (an algebraic group isomorphic to 
some $(\C^*)^{\mu}$), therefore the dual of $\frac{1}{z}\C[\frac{1}{z}]$ 
is a pro-torus; by definition it is the \emph{exponential torus} 
${\bf T}_{exp,int}$ (integral slopes case). Then the tannakian Galois 
group $\pi_{1,form,int}^{\otimes}$ of the category ${\mathcal D}_{form,int}$ 
is isomorphic to the product
of the exponential torus ${\bf T}_{exp,int}$ by the fuchsian group 
$Hom_{gr}(\C^*,\C^*)\times \C$: this is the ``big fundamental group''; 
the ``small fundamental group'' is the product of  the exponential torus 
${\bf T}_{exp,int}$ by the topological fundamental group $\pi_1(D^*,d)$ 
(be careful, the product decompositions \emph{are not canonical}). In the 
general case, without any restriction on the slopes, it is necessary to 
enlarge the universal algebra (replacing the variable $x$ by all its 
ramifications $t^m=x$, $m\in \N^*$). Then there is a \emph{non trivial} 
action of $\gamma$ on the $\Z$-module 
$\displaystyle \bigcup_{m\in \N^*} \frac{1}{x^{1/m}}\C[ \frac{1}{x^{1/m}}]$ 
(by monodromy) and therefore on its dual
${\bf T}_{exp}$, the exponential torus. Then we have in the general case 
\emph{semidirect products}:
$\pi_{1,form}^{\otimes}={\bf T}_{exp}\rtimes (Hom_{gr}(\C^*,\C^*)\times \C)$ 
and $\pi_{1,form}={\bf T}_{exp}\rtimes \pi_1(D^*,d)$. \\

The last step is the study of the category of \emph{meromorphic} 
connections, or equivalently the category ${\mathcal D}_{an}$ of 
$\mathcal D$-modules, where $\mathcal D=\Ka[d/dz]$. This step is 
\emph{very difficult} and involves a lot of delicate and deep analysis. 
Here we shall only describe roughly the results (for more information 
one can read the original papers \cite{MarRam, RamInv1, RamInv2,  RamInv3}, 
and for more details  \cite{Singer-vdPut2}). Heuristically the origin $0$ in $\C$ has an 
``analytic infinitesimal neighborhood'' and an ``algebraic 
infinitesimal neighborhood'', the algebraic neighborhood lying in 
the heart of the analytic neighborhood and being ``very small'' 
(cf. \cite{Deligne86}). The algebraic neighborhood corresponds to
 $\hat {\Z}(1)=\limproj_{n\in \N^*}\mu_n$ ($\mu_n$ is the group 
of complex $n$-th roots of the unity) considered as a \emph{quotient} 
of $Hom_{gr}(\C^*,\C^*)$ (the unipotent component $\C$ corresponds 
to a ``very very small'' neighborhood of $0$ in the heart 
of the algebraic neighborhood).
The fuchsian torus ${\bf T}_f$ corresponds to a ``part very near of 
the algebraic neighborhood''. It remains to understand what happens 
in the ``huge'' region in the analytic neighborhood located between 
the algebraic neighborhood and the exterior, the ``actual world'' $\C^*$: 
one must imagine it as filled by ``points'' (that we shall label $(q,d)$ 
below, $d$ being a direction and $q$ a ``parameter'' of scale). Each point 
will be responsible for a ``monodromy'', the semi-simple part of this 
monodromy will be related to the exponential torus and its unipotent part 
will have an \emph{infinitesimal generator}, which is a Galois derivation: 
we call it an \emph{alien derivation} (and denote it 
$\Der_{q,d}$). It is possible to give a rigorous meaning to 
this heuristic description. There are various approaches, the more 
interesting for the study of $q$-analogues is the tannakian one
(cf. \cite{Deligne86}): one thinks to \emph{fiber functors} as ``points'' 
and to \emph{isomorphisms between fiber functors} as ``paths'' 
(\emph{automorphisms} of fiber functors corresponding to ``loops''). 
Here the paths are made of classical paths (analytic continuation) 
and new paths corresponding to \emph{multisummability} of formal 
(divergent) power series (it is worth noticing that at the algorithmic 
level these two families of paths are in fact very similar: \cite{JNT}). 
Heuristically when you have ``sufficiently many points and loops'', then 
the loops ``fill'' the tannakian Galois group (topologically in Zariski 
sense): this situation will correspond to the ``fundamental group'' 
(the small one). \\

Let us describe now the ``points'' of the ``annulus of the infinitesimal 
neighborhood'' between the algebraic neighborhood and the ``exterior real 
world'' $\C^*$. We shall first give the description and justify it later: 
the points will appear naturally from the analysis of the Stokes phenomena, 
that is from the construction of the paths. We remark that the infinite 
dimensional lattice $\frac{1}{z}\C[\frac{1}{z}]$ is the topological dual 
of its topological dual, the exponential torus
${\bf T}_{exp,int}$. Then each polynomial $q\in \frac{1}{z}\C[\frac{1}{z}]$ 
can be interpreted as a \emph{weight} on the exponential torus: 
if $\tau \in {\bf T}_{exp,int}$, $\tau (e^q)=q(\tau )e^q$,
$q:{\bf T}_{exp,int} \rightarrow \C^*$  is a morphism of pro-algebraic 
groups. The set of directions $d$ issued from the origin is parametrized 
by the unit circle $S^1$ (which we can identify with the boundary of $\C^*$, 
the real blow up of the origin in $\C$ corresponding to $r=0$ in polar 
coordinates $(r,\theta)$). We shall call \emph{degree of} $q$ its degree 
in $1/x$. If $a/x^k$ ($a\in \C^*, k\in \N^*$) is the monomial of highest 
degree of $q$, then it controls the growth or the decay of $e^q$ near 
the origin (except perhaps on the
family of $2k$ ``oscillating lines'': $\Re (a/x^k)=0$, classically named 
Stokes lines or, better..., anti-Stokes lines), we have $k$ open sectors 
of \emph{exponential decay} (of order $k$) of $e^q$ and
$k$ open sectors of \emph{exponential growth} (of order $k$) of $e^q$. 
To each pair 
$(q,d)\in \frac{1}{z}\C[\frac{1}{z}]\times S^1$ such that the direction 
$d$ bisects a sector of decay of $e^q$ we associate a label $(q,d)$: 
the labels will correspond to the points in the ``terra incognita'', 
our mysterious annulus. We introduce on 
$\frac{1}{z}\C[\frac{1}{z}]=\check {{\bf T}}_{exp,int}$ the filtration 
by the degree $k$ (it corresponds to the slope filtration associated 
to the Newton polygon in the formal category). Heuristically, 
if $k=\deg q$, then the corresponding point $(q,d)$ ``belongs'' 
to the direction $d$ and if $k$ is ``big'' this point is far from 
the algebraic neighborhood and near of the exterior world $\C^*$. 
(To each $k\in \N^*$ corresponds a ``slice'' isomorphic to $\C^*$, 
an annulus. If $k>k'$, then the $k'$-annulus is ``surrounded'' by 
the $k$ annulus, and 
``very small'' compared to it  \cite{RamDerniers}). We shall actually 
need points on 
the ``universal covering'' of our annulus. They are labelled by 
the $(q,{\bf d})$, where ${\bf d}$ is a direction above $d$ on the 
Riemann surface of the logarithm. \\

In order to describe the ``paths'', we need the notion of 
\emph{multisummability} (\cite{MaRa}, \cite{Balser}, \cite{Singer-vdPut2}, 
\cite{RamDA}). Let $\hat f \in \Kf$; we shall say that it is \emph{holonomic}
if there exists $D\in {\mathcal D}=\Ka [d/dx]$ such that 
$Df=0$. The set of holonomic power series expansions is a sub-differential 
\emph{algebra} $\mathcal K$ of $\Kf$ (containing $\Ka$) and there is a family 
of summation operators $(S_d^{\pm})_{d\in S^1}$ ($-$ is for ``before $d$'' 
and $+$ is for ``after $d$'' when one turns on $S^1$ in the positive sense):
$S_d^{\pm}: {\mathcal K}\rightarrow {\mathcal O}_d$ (where ${\mathcal O}_d$ 
is the algebra of germs of holomorphic functions on sectors bisected by $d$), 
these operators are \emph{injective homomorphisms of differential algebras}, 
their restriction to $\Ka$ is the classical sum of a convergent power series 
and
$S^{\pm} (\hat f)$ admits $\hat f$ as an \emph{asymptotic expansion};
moreover, for a fixed $\hat f$, the two summations $S_d^+$ and $S_d^-$ 
co\"{\i}ncide, except perhaps for a \emph{finite} set of \emph{singular 
directions}; when $d$ moves between two singular directions the sums 
$S_d^+(\hat f)=S_d^-(\hat f)$ glue together by \emph{analytic continuation}. 
When $d$ crosses a singular line, there is a jump in the sum: this is 
the \emph{Stokes phenomenon}. We consider now the
 differential algebra ${\mathcal U}_{an}=
{\mathcal K}\big\{ (x^{\alpha})_{\alpha\in \C},
(e^q)_{q\in \frac{1}{z}\C[\frac{1}{z}]},\log~x\big\} $, 
it is the \emph{universal differential algebra} associated 
to the family of germs of meromorphic connections.
There are natural extensions of the operators $S_d^{\pm}$ to 
${\mathcal U}_{an}$, but we have to be careful: we must define 
$S_d^{\pm}(\log~x)$ and $S_d^{\pm}(x^{\alpha})$. In order to do 
that we need to choose a branch of the logarithm in a germ of sector 
bisected by $d$ ($x^{\alpha}=e^{\alpha \log~x}$). This corresponds to 
the choice of a direction ${\bf d}$ above $d$ on the Riemann surface of 
the logarithm
(${\bf d}\in (\R,0)$, which is the universal covering of $(S^1,1)$).
In the end, we get a family of summation operators 
$(S_{\bf d}^{\pm})_{{\bf d}\in \R}:
{\mathcal U}_{an}\rightarrow {\mathcal O}_d$;
they are \emph{injective homomorphisms of differential algebras}. \\
Let $\Delta: \frac{dY}{dx}=AY$ be a germ of meromorphic system at 
the origin (integral slopes case). It admits a \emph{formal fundamental 
matrix solution} $\hat F: \frac{d \hat F}{dx}=A\hat F$. The entries 
of $\hat F$ belongs to the universal algebra
${\mathcal U}_{an}$, therefore $F_{\bf d}^+=S_{\bf d}^+(\hat F)$ and 
$F_{\bf d}^-=S_{\bf d}^-(\hat F)$ are germs of \emph{actual} fundamental 
solutions on germs of sectors bisected by $d$. We have
$F_{\bf d}^+=F_{\bf d}^-C_d$, where the \emph{constant} matrix 
$C_d\in Gl_n(\C)$ is a \emph{Stokes matrix} (it is unipotent). 
The map $St_d=(S_{\bf d}^+)^{-1}S_{\bf d}^+$ induces an automorphism 
of the differential algebra $\Ka\{ \hat F\} $, therefore it defines 
an element of the differential Galois group of the system $\Delta$. 
More generally $St_d=(S_{\bf d}^+)^{-1}S_{\bf d}^+$ is an automorphism 
of the  simple differential algebra ${\mathcal U}_{an}$ and defines 
an element of the differential Galois group of this algebra. This element 
is pro-unipotent and we can define a \emph{Galois derivation} 
$\Der_{\bf d}$ of ${\mathcal U}_{an}$ by $St_d=e^{\Der_{\bf d}}$;
by definition, $\Der_{\bf d}$ is the \emph{alien derivation} in 
the direction ${\bf d}$. Now there is a quite \emph{subtle point} 
in our analysis: from the germ of meromorphic system $\Delta: dY/dx=AY$ 
we get a representation $\rho_{\Delta, form}: \pi_{1,form}\rightarrow Gl(V)$ 
and a family of Stokes automorphisms
$(St_{\bf d}(\Delta )\in Gl(V))_{{\bf d}\in \R}$. This last datum is 
equivalent to the knowledge of the corresponding family of alien 
derivations $(\Der_{\bf d}(\Delta )\in End(V))_{{\bf d}\in \R}$. 
There is a natural action of the topological fundamental group on 
the family of alien derivations: 
$\gamma \Der_{\bf d} \gamma^{-1}=\Der_{\gamma ({\bf d})}$ 
($\gamma ({\bf d})$ is a translation of $-2\pi$ of ${\bf d}$), 
therefore it is natural to introduce the semi-direct product 
$\displaystyle exp({\ast}_{{\bf d}\in \R} \C \Der_{\bf d})\rtimes (\gamma )$  
(where
 $\displaystyle \ast_{{\bf d}\in \R} \C \Der_{\bf d})$ 
is the \emph{free Lie-algebra} generated by the symbols 
$\Der_{\bf d}$ and $exp({\ast}_{{\bf d}\in \R} \C \Der_{\bf d})$ 
its exponential group in a ``good sense'') and to observe that 
the connection defined by $\Delta$ ``is'' the representation of 
this group. We could stop here and be happy: why not decide that 
$\displaystyle exp({\ast}_{{\bf d}\in \R} \C \Der_{\bf d})\rtimes (\gamma )$ 
is the fundamental group for the meromorphic category? 
This \emph{does not work}. Of course we have all the knowledge but 
in a bad form: 
to a connection we can associate a representation of our group, 
but conversely there are representations which do not come from 
a connection, the admissible representations are \emph{conditionned}. 
The geometric meaning of the problem is clear: $St_d$ corresponds to 
a loop around a whole bunch of points: all the $(q,{\bf d})$ corresponding 
to all the $q\in \frac{1}{z}\C[\frac{1}{z}]$ admitting $d$ as a line 
of maximal decay for $e^q$ (we shall say in that case that $q$ is 
supported by $d$ and note $(q,{\bf d})\in d$), but a ``good'' fundamental 
group must allow loops around each \emph{individual} point $(q,{\bf d})$. 
It is not difficult to solve the problem; we know \emph{a priori} that our 
representation must contain in some sense the answer, it remains ``only'' 
to extract it. The idea is quite natural: using the exponential torus 
we shall ``vibrate'' the alien derivation $\Der_{\bf d}$ and extracts 
the ``Fourier coefficients'' 
$\Der_{q,{\bf d}}$ (for $q$ supported by $d$). We introduce, in 
the Lie algebra of the differential Galois group $G_{an}$ of 
${\mathcal U}_{an}$,  the family 
$(\tau \Der_{\bf d}\tau^{-1})_{\tau \in {\bf T}_{exp,int}}$ (it is 
a family of Galois derivations), then we consider the ``Fourier expansion''
$\displaystyle \tau \Der_{\bf d}\tau^{-1}=
\sum_{(q,{\bf d})\in d} q(\tau )\Der_{q,{\bf d}}$
(it makes sense because for each connection the sum is finite). 
The coefficients are also in the Lie algebra of the $G_{an}$, 
they are Galois derivations. Now we have won: we consider 
the free Lie algebra 
$\displaystyle Lie~{\mathcal R}=
{\ast}_{{\bf d}\in \R,(q,{\bf d})\in d} \C \Der_{(q,{\bf d})}$
(it is, by definition, the \emph{resurgence algebra}), and the corresponding
exponential group $\mathcal R$ (it makes sense \cite{MaRa}, 
it is by definition the \emph{resurgence group}). We have 
an action of the formal fundamental group on the resurgence Lie algebra:  
$\gamma \Der_{(q,{\bf d})} \gamma^{-1}=\Der_{(q,{\gamma ({\bf d}))}}$,
$\tau \Der_{(q,{\bf d})}\tau^{-1}=
q(\tau ) \Der_{(q,{\bf d})}$ ($\tau \in {\bf T}_{exp,int}$) and
we get a semi-direct product 
$\displaystyle {\mathcal R\rtimes \pi_{1,form,int}=
exp({\ast}_{{\bf d}\in \R,(q,{\bf d})\in d} \C \Der_{(q,{\bf d})})\rtimes 
 ({\bf T}_{exp,int}}\times (\gamma ))$. 
The knowledge of a representation is equivalent to the knowledge 
of its restriction to the formal part and  its ``infinitesimal restriction'' 
to the free Lie algebra. Now the objects of our category (the meromorphic 
connections) correspond to \emph{unconditioned} representations  
(by representation we mean, of course, finite dimensional representation 
whose restriction to the exponential torus is a morphism). We have now 
a fundamental group (the small one), it is the \emph{wild fundamental 
group} (this is in the integral slope case, but with small adaptations 
it is easy to build the wild fundamental group in the general case). 
What about
the big fundamental group (that is the tannakian Galois group)? 
We can easily derive its description from the knowledge of the wild 
fundamental group. The first step is to build some sort of pro-algebraic 
completion of the resurgent Lie algebra $Lie~{\mathcal R}$ 
(cf. \cite{Singer-vdPut2}): if $\rho$ is a representation of our wild 
fundamental group, we can suppose that $V=\C^n$ and that the image of 
the exponential torus is \emph{diagonal}, it follows that the 
corresponding ``infinitesimal restriction''  $\psi =L\rho$ to 
the resurgent Lie algebra satisfies \emph{automatically} the two conditions:
\begin{enumerate}
\item{$\psi (\Der_{(q,{\bf d})})$ is \emph{nilpotent} for every 
$\Der_{(q,{\bf d})}$.}
\item{There are only finitely many $\Der_{(q,{\bf d})}$ such that 
$\psi (\Der_{(q,{\bf d})})\not= 0$.}
\end{enumerate}
By definition, the pro-algebraic completion $(Lie~{\mathcal R})^{alg}$ 
of the free Lie algebra $Lie~{\mathcal R}$ is a projective limit of 
\emph{algebraic}  Lie algebras: 
$(Lie~{\mathcal R})^{alg}=
\displaystyle \limproj _{\psi}~ Lie~{\mathcal R}/\Ker ~ \psi$, 
where the projective limit is taken over all homomorphisms of 
$\C$-algebras $\psi : Lie~{\mathcal R} \rightarrow End(V)$ (where 
$V$ is an arbitrary finite dimensional complex space) satisfying  
conditions (1) and (2). Each algebraic Lie algebra 
$Lie~{\mathcal R}/\Ker ~ \psi$ is the Lie algebra of a connected 
algebraic subgroup of $Gl(V)$. We can consider the projective limit 
of these subgroups, it is a pro-algebraic group (a kind of 
pro-algebraic completion of the resurgent group $\mathcal R$).
We shall call it the resurgent pro-algebraic group and denote it 
${\mathcal R}^{alg}$, its Lie algebra is $(Lie~{\mathcal R})^{alg}$: 
$Lie~{\mathcal R}^{alg}=(Lie~{\mathcal R})^{alg}$. The action of $(\gamma)$ 
on $Lie~{\mathcal R}$ gives an action on $Lie~{\mathcal R}^{alg}$, 
this action can be extended ``by continuity'' to an action of 
$\pi_{1,f}^{\otimes}=Hom_{gr}(\C^*,\C^*)\times \C$, and, using 
the exponential, we get an action of $\pi_{1,f}$ on ${\mathcal R}^{alg}$;
there is also clearly an action of the exponential torus ${\bf T}_{exp,int}$ 
on ${\mathcal R}^{alg}$. Finally, we get a semi-direct product
${\mathcal R}^{alg}\rtimes \pi_{1,form,int}^{\otimes}$, which is 
isomorphic to the tannakian group $\pi_{1,an,int}^{\otimes}$. \\

The aim of this paper is to describe $q$-analogues of the differential
fundamental groups. The construction is also done in three steps: (1)
regular-singular or fuchsian equations, (2) formal or pure equations, (3)
arbitrary equations meromorphic at the origin. We shall limit ourselves to the
integral slopes case (cf. some comments below). We shall recognise the main
actors of the differential case under various disguises. The first two steps
are already well known and the new and difficult part is the last one. The
first step is a little bit more complicated than in the differential case, the
second step is a lot simpler (the exponential torus is replaced by a
\emph{theta torus} and, in the integral slope case, this theta torus is
isomorphic to $\C^*$, therefore radically simpler than the exponential torus).
For the last step, the proofs are less intuitive but  finally the results are
in some sense simpler: one of the essential simplifications is due to the fact
that the $q$-resurgent group is \emph{unipotent} (the differential resurgent
group contains a lot of $Sl_2$ pairs due to the possibility to play with $q$
and $-q$, which exchange the sectors of growth and decay: $e^q,e^{-q}$). \\
We fix $q \in \C$ such that $|q| > 1$ and write 
$q=e^{-2i\pi \tau}$, $\Im~\tau >0$. We begin with the 
\emph{regular singular case}: a germ of meromorphic system at the origin 
$\sigma_q Y=AY$ is regular singular if and only if it is meromophically 
equivalent to a \emph{fuchsian} system $\sigma_q Y=BY$  ($B(0)\in Gl_n(\C)$). 
We call the corresponding category 
${\mathcal E}_f^{(0)}$ the category of \emph{fuchsian modules}, 
its tannakian Galois group is isomorphic to
$Hom_{gr}(\Eq,\C^*) \times \C$, where $\Eq=\C^*/q^{\Z}$ is the elliptic 
curve associated to $q$ (cf. \cite{JSGAL} 2.2.2). There exists also a 
``small group'' Zariski dense in the tannakian Galois group, and one can guess 
it using a $q$-analogy: the image of $\Z$ in $Hom_{gr}(\C^*,\C^*)$ is 
the subgroup of group homomorphisms which are algebraic group homomorphisms, 
therefore it is natural to consider the subgroup $\Pi$ of the elements 
of $Hom_{gr}(\Eq,\C^*)$ which are \emph{continuous}. We use 
the decomposition $\C^*= \U\times q^{\R}$ ($\U \subset \C^{*}$ 
is the unit circle) and we denote 
$\gamma_1, \gamma_2 \in Hom_{gr}(\Eq,\C^*)$ the continuous group 
homomorphisms defined respectively by $uq^y \mapsto u$ and 
$uq^y \mapsto e^{2i\pi y}$. Then $\Pi$ is generated by $\gamma_1$ 
and $\gamma_2$ and is Zariski-dense in $Hom_{gr}(\C^*,\C^*)$,  the ``fundamental group"  of the category ${\mathcal E}_f^{(0)}$ (the local fundamental group) $\pi_{1,q,f}$  is by definition the subgroup of $Hom_{gr}(E_q,\C^*\times \C)$ whose semi-simple component is generated by $\gamma_1$ and $\gamma_2$ and whose unipotent component is $\Z$ (cf. \cite{JSGAL} 2.2.2).\\

The next step is the study of the category $\EE_{form}$ of 
\emph{formal} $q$-difference modules. We shall limit ourselves to the 
integral slope case: the category $\EE_{form,int}$ (or equivalently 
of the category $\EE_{p,1}^{(0)}$ of pure meromorphic modules 
with integral slopes, cf. below). It is a neutral tannakian category. 
As in the differential case, in order to compute the corresponding 
``fundamental groups'', it is necessary to understand the formal 
classification of $q$-difference equations of order one: two such 
equations $\sigma_qy-\hat ay=0$ and 
$\sigma_qy-\hat by=0$ ($\hat a,\hat b\in \Kf^*$) 
are formally equivalent if and only if $a^{-1}b\in \sigma_{q,\log}~\Kf$, 
where $\sigma_{q,\log}\hat f=\sigma_q(\hat f)/\hat f$. Then the order one 
equations are classified by the abelian group 
$\C^*/q^{\Z}\times (z^m)_{m\in \Z}\simeq \Eq\times \Z$  
($\Eq$ correspond to the fuchsian equations, $(z^m)_{m\in \Z}$ 
to irregular equations). The ``basic'' irregular equation is
$\sigma_q y-zy=0$, it admits the Jacobi theta function $\theta_q$ as 
a solution (cf. below) and its $q$-difference Galois group is isomorphic 
to $\C^*$. Then one can prove that the tannakian Galois group 
$G_{form,int}$ of the category $\EE_{form,int}$ is isomorphic 
to the topological dual group of $ \Eq\times \Z$ 
(where $\Eq$ is interpreted as the inductive limit of its finitely 
generated subgroups), that is to
$\C^*\times (Hom_{gr}(\Eq,\C^*)\times \C)$: $\C^*$ is by definition the
\emph{theta torus}, it is the $q$-analogue of the exponential torus. (The
tannakian Galois group $G_{p,1}^{(0)}$ of the category 
$\EE_{p,1}^{(0)}$ of pure meromorphic modules with integral slope is
isomorphic to $G_{form,int}$.) We do not know what will happen in the non
integral slope case.\\

The last step and the main purpose of this paper is the study of  
the category $\EE_1^{(0)}$ of $q$-difference modules whose Newton 
polygon admits only integral slopes. It is a neutral tannakian category, 
we shall prove that there exists a \emph{semi-direct} decomposition 
of its tannakian Galois group $G_{1}^{(0)} = \St \rtimes G_{p,1}^{(0)}$, 
where $\St$ is a unipotent pro-algebraic group, and we shall describe 
the Lie algebra $\st$ of $\St$: like in the differential case this Lie 
algebra is a ``pro-algebraic completion'' of a \emph{free} complex 
Lie algebra generated by a family of $q$-alien derivations: 
$\big(\Derdc\big) _{\delta \in \N^{*}, \bar a \in \Eq}$. \\
These $q$-alien derivations are indexed by labels $(\delta,{\bf a})$ which 
are the $q$-analogs of the labels $(q,{\bf d})$ of the differential case:
$\delta$ is a weight on the $\theta$-torus $\C^*$ (that is, an element of 
the topological dual group $\Z$; actually, only the $\delta > 0$ have a
non trivial action, so that we harmlessly take $\delta \in \N^{*}$), 
and $\bf a$  is a pair formed by 
$a \in E_q$ (\emph{i.e.} a $q$-direction, representing a germ of $q$-spiral 
at the origin) and an element $\xi$ of the $q$-local fundamental group 
$\pi_{1,q,f}$. In order to define the $q$-alien derivations, we will 
use, as in the differential case, some summability tools (here, an algebraic 
version of the $q$-multisummability due to second author), but the approach 
will be different: we will no longer use solutions but replace them by 
fiber functors ($\otimes$-functors). We will deal with meromorphic 
families of Lie-like automorphisms of fiber functors (the variable being 
the $q$-direction of summability) and extract their singularities by 
a residue process, this will give birth to $q$-alien derivations 
$\dot \Delta_a$. \\

In this paper, we shall compute the $q$-alien derivations in the one-level 
case using a $q$-Borel transform (of some convenient order). This relates 
alien derivations to the irregularity invariants introduced in \cite{RSZ} 
and proves that, in this case, $q$-alien derivations are \emph{a complete 
set of irregularity invariants}. We shall extend these results to the general 
case in a forthcoming paper \cite{RS2}. The principle is similar but it 
is necessary to introduce a double family of categories ``interpolating'' 
between $\EE_1^{(0)}$ and respectively ${\mathcal E}_{form}$ 
and $\EE_{p,1}^{(0)}$, in relation with slopes and $q$-Gevrey estimates. 
With these tools, it is possible to prove that, in the general 
case also, $q$-alien derivations are a complete set of irregularity 
invariants and that the $q$-resurgence group is Zariski dense in $\St$. 
(The reader can check as an exercise that, in the general case, for an 
isoformal family of meromorphic $q$-difference modules $M$, the dimensions 
of the $\C$-vector space of the irregular invariants of \cite{RSZ} and of 
the $\C$-vector space generated by the ``acting" $q$-alien derivations 
are equal: they are equal to the area of the ``closed Newton polygon" of $M$).



\section{Prerequisites (mostly from \cite{JSGAL}, \cite{JSFIL}, \cite{RSZ}
and \cite{JSSTO})}


\subsection{General facts}
\label{subsection:generalfacts}

\paragraph{Notations, general conventions.}
We fix $q \in \C$ such that $|q| > 1$. We then define the 
automorphism $\sq$ on various rings, fields or spaces of functions 
by putting $\sq f(z) = f(qz)$. This holds in particular for 
the ring $\Ra$ of convergent power series and its field 
of fractions $\Ka$, the ring $\Rf$ of formal power series 
and its field of fractions $\Kf$, the ring $\Rwg$ of holomorphic 
germs and the field $\Kwg$ of meromorphic germs in the punctured 
neighborhood of $0$, the ring $\Rw$ of holomorphic functions
and the field $\Kw$ of meromorphic functions on $\C^{*}$; 
this also holds for all modules or spaces of vectors
or matrices over these rings and fields. For any such ring 
(resp. field) $R$, the $\sq$-invariants elements
make up the subring (resp. subfield) $R^{\sq}$ of constants. 
The field of constants of $\M(\C^{*},0)$ and that of
$\M(\C^{*})$ can be identified with a field of elliptic functions,
the field $\M(\Eq)$ of meromorphic functions over the complex torus
$\Eq = \C^{*}/q^{\Z}$. \\
We shall write $\overline{a} = \pi(a) \in \Eq$ for the image of 
$a \in \C^{*}$ by the natural projection $\pi: \C^{*} \rightarrow \Eq$, 
and $[a;q] = a q^{\Z} = \pi^{-1}(\overline{a}) \subset \C^{*}$ 
for the preimage of $\overline{a}$ in $\C^{*}$, a discrete $q$-spiral. 
These notations extend to subsets $A \subset \C^{*}$: 
$\overline{A} = \pi(A) \subset \Eq$ and 
$[A;q] = A q^{\Z} = \pi^{-1}(\overline{A}) \subset \C^{*}$.

\paragraph{Categories.}
Let $K$ denote any one of the forementioned fields of functions.
Then, we write $\D = K\l<\sigma,\sigma^{-1}\r>$ for the \"{O}re algebra
of non commutative Laurent polynomials characterized by the relation
$\sigma . f = \sq(f) . \sigma$. We now define the category of 
$q$-difference modules in three clearly equivalent ways:
\begin{eqnarray*}
\DM & = & \{(E,\Phi) \;/\; E \text{~a~} K \text{-vector space of finite rank~},
                        \Phi : E \rightarrow E \; 
                               \text{~a~} \sq \text{-linear} \text{~map}\} \\
    & = & \{(K^{n}, \Phi_{A}) \;/\; A \in GL_{n}(K) ,
                                    \Phi_{A}(X) = A^{-1} \sq X \} \\
    & = & \{\text{~finite length left~} \D \text{-modules}\}.
\end{eqnarray*}
For instance, a morphism from $M_{A} = (K^{n}, \Phi_{A})$ to 
$M_{B} = (K^{n}, \Phi_{B})$, where $A \in GL_{n}(K)$ and $B \in GL_{p}(K)$, 
is a $F \in M_{p,n}(K)$ such that $(\sq) F A = B F$. Then, $\DM$ is 
a $\C$-linear abelian rigid tensor category, hence a tannakian category. 
Moreover, all objects in $\DM$ have the form $\D / \D P$. In the case
of $K = \Ka$, the category $\DM$ will be written $\EE^{(0)}$
(for ``equations near $0$'').

\paragraph{Vector bundles and fiber functors.}
To any module $M_{A}$ in $\EE^{(0)}$, one can associate a holomorphic
vector bundle $\F_{A}$ over $\Eq$:
$$
\F_{A} = \dfrac{(\C^{*},0) \times \C^{n}}{(z,X) \sim (q z,A(z) X)} 
\rightarrow \dfrac{(\C^{*},0)}{z \sim q z} = \Eq.
$$
This is the usual construction from equivariant bundles except that
the germ $(\C^{*},0)$ is only endowed with the action of the semigroup
$q^{-\N}$ instead of a group; correspondingly, the projection map
is not a covering. The pullback $\tilde{\F_{A}} = \pi^{*}(\F_{A})$
over the open Riemann surface $\C^{*}$ is the trivial bundle 
$\C^{*} \times \C^{n}$, but with an equivariant action
by $q^{\Z}$. The $\O_{\Eq}$-module of sections of $\F_{A}$ (also
written $\F_{A}$) is the sheaf over $\Eq$ defined by:
$\F_{A}(V) = 
\{\text{solutions of~} \sq X = A X \text{~holomorphic over~} \pi^{-1}(V)\}$.
>From these two descriptions, the following is immediate:
\begin{prop}
This gives an exact faithful $\otimes$-functor $M_{A} \leadsto \F_{A}$
from $\EE^{(0)}$ to the category $Fib(\Eq)$ of holomorphic vector bundles.
Taking the fiber of $\tilde{\F_{A}}$ at $a \in \C^{*}$ yields a
fiber functor $\omega^{(0)}_{a}$ on $\EE^{(0)}$ over $\C$.
\end{prop}

\paragraph{Newton polygon.}
Any $q$-difference module $M$ over $\Ka$ or $\Kf$, can be given a
Newton polygon $N(M)$ at $0$, or, equivalently, a Newton function
$r_{M}$ sending the slope 
\footnote{It should be noted that the slopes defined and used in 
the present paper are the \emph{opposites} of the slopes defined 
in previous papers.} 
$\mu \in S(M) \subset \Q$ to its multiplicity $r_{M}(\mu) \in \N^{*}$ 
(and the $\mu$ out of the support $S(M)$ to $0$). \\ 
For instance, the $q$-difference operator 
$L = q z \sigma^{2} - (1 + z) \sigma + 1 \in \D$ gives rise to 
the $q$-difference equation $q z \sq^{2} f - (1+z) \sq f+ f = 0$, 
of which the so-called Tschakaloff series 
$\sum_{n \geq 0}  q^{n(n-1)/2} z^{n}$ is a solution (it is a
natural $q$-analogue of the Euler series). By vectorisation, this
equation gives rise to the system $\sq X = A X$, where 
$A = \begin{pmatrix} z^{-1} & z^{-1} \\ 0 & 1 \end{pmatrix}$,
and to the module $M = M_{A}$. The latter is isomorphic to
$\D/\D \hat{L}$, where 
$\hat{L} = \sigma^{2} - (z + 1) \sigma + z = (\sigma - z)(\sigma - 1)$
is the dual operator of $L$. We respectively attach to $\sigma - z$ 
and $\sigma - 1$ the slopes $-1$ and $0$ and take $S(M) = \{-1,0\}$, 
with multiplicities $r_{M}(-1) = r_{M}(0) = 1$.


\subsection{Filtration by the slopes}
\label{subsection:filtration}

\paragraph{The filtration and the associated graded module.}
The module $M$ is said to be pure isoclinic of slope $\mu$ if 
$S(M) = \{\mu\}$ and fuchsian if moreover $\mu = 0$.
Direct sums of pure isoclinic modules are called pure 
modules~\footnote{It should be noted that we call in the present paper 
a pure isoclinic (resp. pure) module what was called a pure 
(resp. tamely irregular) module in previous papers.}: they are 
irregular objects without wild monodromy, as follows from \cite{RSZ} 
and \cite{JSSTO}. 
The tannakian subcategory of $\EE^{(0)}$ made up of pure modules
is called $\EE^{(0)}_{p}$. Modules with integral slopes also form
tannakian subcategories, which we write $\EE^{(0)}_{1}$ and 
$\EE^{(0)}_{p,1}$. \\
It was proved in \cite{JSGAL} that the category $\EE^{(0)}_{f}$
of fuchsian modules is equivalent to the category of \emph{flat} holomorphic
vector bundles over $\Eq$ and that its Galois group $G_{f}^{(0)}$
is isomorphic to $Hom_{gr}(\C^{*}/q^{\Z},\C^{*}) \times \C$
(here, as in the introduction, $Hom_{gr}$ means ``morphisms of \emph{abstract} groups'').
Since objects of $\EE^{(0)}_{p,1}$ are essentially $\Z$-graded
objects with fuchsian components, the Galois group of $\EE^{(0)}_{p,1}$
is $G_{p,1}^{(0)} = \C^{*} \times G_{f}^{(0)}$.
\begin{thm}\cite{JSFIL}
\label{theorem:filtration}
Let the letter $K$ stand for the field $\Ka$ (convergent case)
or the field $\Kf$ (formal case). In any case, any object $M$ of $\DM$ 
admits a unique filtration $(F_{\leq \mu}(M))_{\mu \in \Q}$ by subobjects
such that each $F_{(\mu)}(M) = \frac{F_{\leq \mu}(M)}{F_{< \mu}(M)}$ 
is pure of slope $\mu$ (thus of rank $r_{M}(\mu)$). The $F_{(\mu)}$ are
endofunctors of $\DM$ and $\gr = \bigoplus F_{(\mu)}$ is a faithful exact 
$\C$-linear $\otimes$-compatible functor and a retraction of the inclusion 
of $\EE_{p}^{(0)}$ into $\EE^{(0)}$. In particular, the functor $\gr$ 
retracts $\EE_{1}^{(0)}$ to $\EE_{p,1}^{(0)}$.In the formal case, 
$gr$ is isomorphic to the identity functor.
\end{thm}
\begin{cor}
For each $a \in \C^{*}$, the functor 
$\hat{\omega}_{a}^{(0)} = \omega_{a}^{(0)} \circ \gr$ is a fiber functor.
\end{cor}
We shall consistently select an arbitrary basepoint $a \in \C^{*}$ 
and identify the Galois group $G^{(0)}$  as 
$Aut^{\otimes}\bigl(\hat{\omega}_{a}^{(0)}\bigr)$.
\begin{cor}
The Galois group $G^{(0)}$ of $\EE^{(0)}$ is the semi-direct product
$\St \rtimes G_{p}^{(0)}$ of the Galois group $G_{p}^{(0)}$ of 
$\EE_{p}^{(0)}$ by a prounipotent group, the Stokes group $\St$.
\end{cor}
\emph{\textbf{From now on, we only consider modules with integral slopes}}. \\
Further studies would have to be based on the work \cite{vdPR} by
van der Put and Reversat. 

\paragraph{Description in matrix terms.}
We now introduce notational conventions which will be used all along 
this paper for a module $M$ in $\EE_{1}^{(0)}$ and its associated graded
module $M_{0} = gr(M)$, an object of $\EE_{p,1}^{(0)}$. The module $M$ 
may be given the shape 
$M = (\Ka^{n},\Phi_{A})$, with:
\begin{equation} 
\label{equation:forme-canonique}
A = A_{U} \underset{def}{=}
\begin{pmatrix}
z^{\mu_{1}} A_{1}  & \ldots & \ldots & \ldots & \ldots \\
\ldots & \ldots & \ldots  & U_{i,j} & \ldots \\
0      & \ldots & \ldots   & \ldots & \ldots \\
\ldots & 0 & \ldots  & \ldots & \ldots \\
0      & \ldots & 0       & \ldots & z^{\mu_{k}} A_{k}    
\end{pmatrix},
\end{equation}
where the slopes $\mu_{1} < \cdots < \mu_{k}$ are integers,
$r_{i} \in \N^{*}$, $A_{i} \in GL_{r_{i}}(\C)$ ($i = 1,\ldots,k$) and
$$
U = (U_{i,j})_{1 \leq i < j \leq k} \in 
\underset{1 \leq i < j \leq k}{\prod} \Ma_{r_{i},r_{j}}(\Ka).
$$
The associated graded module is then a direct sum
$M_{0} = P_{1} \oplus \cdots \oplus P_{k}$, where,
for $1 \leq i < j \leq k$, the module $P_{i}$ is pure 
of rank $r_{i}$ and slope $\mu_{i}$ and can be put into
the form $P_{i} = (\Ka^{r_{i}},\Phi_{z^{\mu_{i} A_{i}}})$.
Therefore, one has $M_{0} = (\Ka^{n},\Phi_{A_{0}})$, where
the matrix $A_{0}$ is block-diagonal (it is the same as $A_{U}$,
with all $U_{i,j} = 0$). \\
We write $\G \subset GL_{n}$ for the algebraic subgroup 
and $\g$ for its Lie algebra, made up of matrices of the form:
\begin{equation} \label{equation:automorphisme}
F = 
\begin{pmatrix}
I_{r_{1}} & \ldots & \ldots & \ldots & \ldots \\
\ldots & \ldots & \ldots  & F_{i,j} & \ldots \\
0      & \ldots & \ldots   & \ldots & \ldots \\
\ldots & 0 & \ldots  & \ldots & \ldots \\
0      & \ldots & 0       & \ldots & I_{r_{k}}     
\end{pmatrix}
\quad \text{~and~} \quad
f = 
\begin{pmatrix}
0_{r_{1}} & \ldots & \ldots & \ldots & \ldots \\
\ldots & \ldots & \ldots  & F_{i,j} & \ldots \\
0      & \ldots & \ldots   & \ldots & \ldots \\
\ldots & 0 & \ldots  & \ldots & \ldots \\
0      & \ldots & 0       & \ldots & 0_{r_{k}}     
\end{pmatrix}.
\end{equation}
For $F$ in $\G$, we shall write $F[A] = \l(\sq F\r) A F^{-1}$
the result of the gauge transformation $F$ on the matrix $A$. 
Theorem \ref{theorem:filtration} entails:
$$
\forall (U_{i,j})_{1 \leq i < j \leq k} \in 
\prod_{1 \leq i < j \leq k} \Ma_{r_{i},r_{j}}(\Ka) \;,\;
\exists ! \hat{F} \in \G(\Kf) \;:\;
\hat{F} [A_{0}] = A_{U}.
$$
This $\hat{F}$ will be written $\hat{F}_{A}$ (where $A = A_{U}$).
The blocks $\hat{F}_{i,j}$ are recursively computed as follows.
For $j < i$, $\hat{F}_{i,j} = 0$. For $j = i$, $\hat{F}_{i,j} = I_{r_{i}}$.
Then, for $j > i$, one must solve the non homogeneous first order
equation:
\begin{equation}
\label{equation:calculrecursifhat(F)}
\sq \hat{F}_{i,j} z^{\mu_{j}} A_{j} - z^{\mu_{i}} A_{i} \hat{F}_{i,j} =
\sum_{i < k < j} U_{i,k} \hat{F}_{k,j} + U_{i,j}.
\end{equation}

\paragraph{Description of the Stokes group.}
To go further, we choose to fix an arbitrary basepoint $a \in \C^{*}$
(see the corollary to theorem \ref{theorem:filtration}) and we identify 
the Galois groups accordingly:
$$
G_{1}^{(0)} \underset{def}{=} Gal(\EE_{1}^{(0)}) = 
Aut^{\otimes}(\hat{\omega}^{(0)}_{a}).
$$
Recall from the quoted papers the action of the pure component
$G_{p,1}^{(0)} = \C^{*} \times Hom_{gr}(\C^{*}/q^{\Z},\C^{*}) \times \C$.
We keep the notations above. For any $A$ with graded part $A_{0}$,
an element $(\alpha,\gamma,\lambda) \in G_{p,1}^{(0)}$ yields the
automorphism of 
$\hat{\omega}^{(0)}_{a}(A) = \omega^{(0)}_{a}(A_{0}) = \C^{n}$
given by the matrix:
$$
\begin{pmatrix}
\alpha^{\mu_{1}} \gamma(A_{s,1}) A_{u,1}^{\lambda}  
                                        & \ldots & \ldots & \ldots & \ldots \\
\ldots & \ldots & \ldots  & 0 & \ldots \\
0      & \ldots & \ldots   & \ldots & \ldots \\
\ldots & 0 & \ldots  & \ldots & \ldots \\
0      & \ldots & 0       & \ldots & 
                           \alpha^{\mu_{k}} \gamma(A_{s,k}) A_{u,k}^{\lambda}
\end{pmatrix},
$$
where we have written $A_{i} = A_{s,i} A_{u,i}$ the multiplicative
Dunford decomposition (into a semi-simple and a unipotent factor that
commute) and $\gamma$ acts on a semi-simple matrix through its eigenvalues. 
These matrices generate the group $G_{p,1}^{(0)}(A) \subset GL_{n}(\C)$. \\
We then have a semi-direct decomposition: 
$G_{1}^{(0)} = \St \rtimes G_{p,1}^{(0)}$, where the Stokes group $\St$ 
is the kernel of the morphism $G_{1}^{(0)} \rightarrow G_{p,1}^{(0)}$.
The group $\St(A)$ is an algebraic subgroup of $\G(\C)$. The above
matrix of $G_{p,1}^{(0)}(A)$ acts by conjugation on the matrix
described by (\ref{equation:automorphisme}): the $F_{i,j}$ block
is sent to 
$\alpha^{\mu_{i} - \mu_{j}} \gamma(A_{s,i}) A_{u,i}^{\lambda}
F_{i,j} \bigl(\gamma(A_{s,j}) A_{u,j}^{\lambda}\bigr)^{-1}$.
In particular, the group $\C^{*}$ acts on the ``level $\delta$''
upper diagonal $\mu_{j} - \mu_{i} = \delta$ (where $\delta \in \N$)
by multiplication by $\alpha^{-d}$. The group $\St(A)$ is
\emph{filtered} by the \emph{normal} subgroups $\St_{\delta}(A)$ defined by:
$\mu_{j} - \mu_{i} \geq \delta$ (meaning that all blocks such that
$0 < \mu_{j} - \mu_{i} < \delta$ vanish). \\
Likewise, the Lie algebra $\st(A) = Lie\bigl(\St(A)\bigr)$,
which is a subalgebra of $\g(\C)$, admits an adjoint action
described by the same formulas (this is because 
$\log P F P^{-1} = P \log F P^{-1}$). The algebra $\st(A)$ 
is \emph{graded} by its ``level $\delta$'' upper diagonals $\st_{\delta}(A)$,
defined by $\mu_{j} - \mu_{i} = \delta$. As noted in \cite{JSSTO},
the algebra $\st_{\delta}(A)$ can be identified with the (group) kernel 
of the central extension 
$\St(A)/\St_{\delta+1}(A) \rightarrow \St(A)/\St_{\delta}(A)$.


\subsection{Stokes operators}
\label{subsection:Stokesoperators}

\paragraph{Algebraic summation.}
The following computations are extracted from \cite{JSSTO}.
We need the following theta function of Jacobi:
$\th_{q}(z) = \sum_{n \in \Z} q^{-n(n+1)/2} z^{n}$.
It is holomorphic in $\C^{*}$ with simple zeroes, all
located on the discrete $q$-spiral $[-1;q]$. It satisfies the
functional equation: $\sq \th_{q} = z \th_{q}$. We then define
$\th_{q,c}(z) = \th_{q}(z/c)$ (for $c \in \C^{*}$); it is holomorphic 
in $\C^{*}$ with simple zeroes, all located on the discrete 
$q$-spiral $[-c;q]$ and satisfies the functional equation: 
$\sq \th_{q,c} = \frac{z}{c} \th_{q,c}$. \\
For a given formal class described by $A_{0}$, $\mu_{1},\ldots,\mu_{k}$
and $r_{1},\ldots,r_{k}$ as above, and for any $c \in \C^{*}$,
we introduce the matrix $T_{c,A_{0}} \in \Ma_{n}(\Kw)$ which
is block-diagonal with blocks $\theta_{c}^{-\mu_{i}} I_{r_{i}}$.
Moreover, we shall assume the following normalisation due to Birkhoff 
and Guenther (see \cite{RSZ}, \cite{JSSTO}):
$$
\forall i < j \;,\; \text{~all coefficients of $U_{i,j}$ belong to~}
\sum_{\mu_{i} \leq d < \mu_{j}} \C z^{d}.
$$
Then, putting $A'_{i} = c^{\mu_{i}} A_{i} \in GL_{r_{i}}(\C)$ and
$U'_{i,j} = (z/c)^{-\mu_{i}} \theta_{c}^{\mu_{j} -\mu_{i}} U_{i,j}
\in \Ma_{r_{i},r_{j}}(\Rw)$, we have:
$$
A'_{U'} \underset{def}{=}
T_{c,A_{0}}[A_{U}] = \begin{pmatrix}
A'_{1}  & \ldots & \ldots & \ldots & \ldots \\
\ldots & \ldots & \ldots  & U'_{i,j} & \ldots \\
0      & \ldots & \ldots   & \ldots & \ldots \\
\ldots & 0 & \ldots  & \ldots & \ldots \\
0      & \ldots & 0       & \ldots & A'_{k}    
\end{pmatrix}.
$$
Now, if the images in $\Eq$ of the
spectra $Sp(A'_{i})$ are pairwise disjoint, there is a unique 
$F' \in \G(\Rw)$ such that $F'[A'_{0}] = A'_{U'}$. Its coefficients 
are recursively defined by the equations:
$$
\sq F'_{i,j} A'_{j} - A'_{i} F'_{i,j} =
\sum_{i < k < j} U'_{i,k} F'_{k,j} + U'_{i,j}.
$$
The unique solution of this equation in $\Rw$ is obtained by taking 
the Laurent series:
$$
F'_{i,j} = \sum_{p \in \Z} \Phi_{q^{p} A'_{j},A'_{i}}^{-1}(V_{p}) \, z^{p},
\quad \bigl(\sum_{p \in \Z} V_{p} \, z^{p} = 
\sum_{i < k < j} U'_{i,k} F'_{k,j} + U'_{i,j}\bigr),
$$
where one writes $\Phi_{B,C}(M) = M B - C M$ (that map is one to one
if and only if $Sp(B) \cap Sp(C) = \emptyset$). Note for further use 
that the condition we have to impose on the spectra is the following:
\begin{equation}
\label{equation:nonresonance}
\forall i < j \;,\; 
q^{\Z} c^{\mu_{i}} Sp(A_{i}) \cap q^{\Z} c^{\mu_{j}} Sp(A_{j}) = \emptyset.
\end{equation}
This is equivalent to requiring that $\overline{c} \not\in \Sigma_{A_{0}}$,
where $\Sigma_{A_{0}}$ is some explicit finite subset of $\Eq$. \\
>From the equalities $A'_{U'} = T_{c,A_{0}}[A_{U}]$,
$A'_{0} = T_{c,A_{0}}[A_{0}]$ and $F'[A'_{0}] = A'_{U'}$,
we get at last $F[A_{0}] = A_{U}$, where
$F = T_{c,A_{0}}^{-1} F' T_{c,A_{0}}$ can be easily
computed: it belongs to $\G(\Kw)$ and 
$F_{i,j} = \theta_{c}^{\mu_{i} - \mu_{j}} F'_{i,j}$. The condition
that the $F'_{i,j}$ are holomorphic over $\C^{*}$ is equivalent to the
following condition:
\begin{equation}
\label{equation:conditionpolaire}
\forall i < j \;,\; F_{i,j} \text{~has poles only on~} [-c;q], 
\text{~and with multiplicities~} \leq \mu_{j} - \mu_{i}.
\end{equation}
Then, we get the following conclusion: there is a \emph{unique} $F \in \G(\Kw)$ 
such that condition (\ref{equation:conditionpolaire}) holds and 
$F[A_{0}] = A_{U}$. Note that the condition depends on 
$\overline{c} \in \Eq$ rather than $c$. To summarize the discussion:
\begin{prop}
For every $c \in \C^{*}$ satisfying condition (\ref{equation:nonresonance})
(\emph{i.e.}, $\overline{c} \not\in \Sigma_{A_{0}}$), there is a unique
$F \in \G(\Kw)$ satisfying condition (\ref{equation:conditionpolaire}) 
and such that $F[A_{0}] = A_{U}$. We consider this $F$ as obtained by 
\emph{summation of $\hat{F}_{A}$ in the direction $\overline{c} \in \Eq$} 
and, accordingly, write it $S_{\overline{c}} \hat{F}_{A}$.
\end{prop}
If we now choose two \emph{$q$-directions of summation}
$\overline{c},\overline{d} \in \Eq$, the ambiguity of summation
is expressed by:
\begin{equation}
\label{equation:defoperateurStokes}
S_{\overline{c},\overline{d}} \hat{F}_{A} \underset{def}{=}
(S_{\overline{c}} \hat{F}_{A})^{-1} S_{\overline{d}} \hat{F}_{A}.
\end{equation}
This is a meromorphic automorphism of $A_{0}$. As explained in \cite{RSZ} 
and \cite{JSSTO}, it is a \emph{Stokes operator}.

\paragraph{Case of one level.}
\label{subsubsection:caseofonelevel}
For further use, we now specialize some of the previous results to
the case of two (integral) slopes $\mu < \nu$, and only one ``level''
$\delta = \nu - \mu \in \N^{*}$. For simplicity, we write our matrices:
\begin{equation}
\label{equation:onelevelmatrices}
M_{0} = \begin{pmatrix} z^{\mu} A & 0 \\ 0 & z^{\nu} B \end{pmatrix}
\quad \text{and} \quad
M = \begin{pmatrix} z^{\mu} A & z^{\mu} U B \\ 0 & z^{\nu} B \end{pmatrix},
\end{equation}
where 
$A \in GLr(\C)$, $B \in GL_{s}(\C)$ and $U \in \Ma_{r,s}(\C_{\delta-1}[z])$
(polynomials with degree $< \delta$). It is clear that the upper right
block can indeed be written in such a way. 
Then the unique element of $\G(\Kf)$ which sends $M_{0}$ to $M$ is 
the matrix $\begin{pmatrix} I_{r} & F \\ 0 & I_{s} \end{pmatrix}$, 
where $F$ is the unique element of $\Ma_{r,s}(\Kf)$ such that:
\begin{equation}
\label{equation:onelevelequation1}
z^{\delta} \sq F - \Lambda(F) = U,
\end{equation}
Here, we have written $\Lambda(F) = A F B^{-1}$ (thus, an endomorphism
of $\Ma_{r,s}(\C)$ and similar spaces). The formal solution $F$
can be computed by identification of coefficients, \emph{i.e.} by solving:
\begin{equation}
\label{equation:onelevelequation2}
\forall n \in \Z \;,\; 
q^{n-\delta} F_{n-\delta} - A F_{n} B^{-1} = U_{n}.
\end{equation}
Similarly, the unique element of $\G(\Kw)$ such that condition 
(\ref{equation:conditionpolaire}) holds which sends $M_{0}$ to $M$ is 
the matrix $\begin{pmatrix} I_{r} & F \\ 0 & I_{s} \end{pmatrix}$, 
where $F$ is the unique of $\Ma_{r,s}(\Kw)$ with poles only on $[-c;q]$ 
and with multiplicities $\leq d$ which is solution of equation
(\ref{equation:onelevelequation1}). This is solved by putting
$F = \theta_{c}^{-\delta} G$, so that $G$ is a solution holomorphic
on $\C^{*}$ of the following equation:
\begin{equation}
\label{equation:onelevelequation3}
c^{\delta} \sq G - \Lambda(G) = V, \text{~where~} 
G = \theta_{c}^{\delta} F \text{~and~} V = \theta_{c}^{\delta} U.
\end{equation}
This can be solved by identification of coefficients of the corresponding
Laurent series, \emph{i.e.} by solving:
\begin{equation}
\label{equation:onelevelequation4}
\forall n \in \Z \;,\; 
c^{\delta} q^{n} G_{n} - A G_{n} B^{-1} = V_{n}.
\end{equation}
This is possible if $c^{\delta} q^{\Z} \cap Sp(A)/Sp(B) = \emptyset$,
which is precisely condition (\ref{equation:nonresonance}) specialized
to the present setting.
Then we can take $G_{n} = (c^{\delta} q^{n} \Id - \Lambda)^{-1} V_{n}$.



\section{Stokes operators and alien derivations}


\subsection{Stokes operators are galois}

We take on the notations of section \ref{subsection:filtration}
and consider moreover another object $B$, to which we apply
similar notations: graded $B_{0}$, diagonal blocks $B_{j}$
corresponding to slopes $\nu_{j}$ with multiplicities $s_{j}$, etc.
\begin{lemma}
\label{lemma:tensfonctrestreint}
(i) Assume the following condition:
\begin{equation}
\label{equation:nonresonancemixtetensorielle}
\forall i < i', j < j' \;,\; 
q^{\Z} c^{\mu_{i} + \nu_{j}} Sp(A_{i}) Sp(B_{j}) \cap 
q^{\Z} c^{\mu_{i'} + \nu_{j'}} Sp(A_{i'}) Sp(B_{j'}) = \emptyset.
\end{equation}
Then:
\begin{equation}
\label{equation:sommationtensorielle}
S_{\overline{c}} \hat{F}_{A \otimes B} = 
(S_{\overline{c}} \hat{F}_{A}) \otimes (S_{\overline{c}} \hat{F}_{B}).
\end{equation}
(ii) Assume the following condition:
\begin{equation}
\label{equation:nonresonancemixtefonctorielle}
\forall i, j \text{~such that~} \mu_{i} < \nu_{j} \;,\; 
q^{\Z} c^{\mu_{i}} Sp(A_{i}) \cap q^{\Z} c^{\nu_{j}} Sp(B_{j}) = \emptyset.
\end{equation}
Then, for any morphism $F : A \rightarrow B$, writing $F_{0} = \gr F$, we have:
\begin{equation}
\label{equation:sommationfonctorielle}
F \; S_{\overline{c}} \hat{F}_{A} = S_{\overline{c}} \hat{F}_{B} \; F_{0}.
\end{equation}
\end{lemma}
\Pr
(i) From elementary properties of the tensor product, we draw that
the diagonal blocks of $A \otimes B$ are the 
$z^{\mu_{i} + \nu_{j}} A_{i} \otimes B_{j}$ and that
$Sp(A_{i} \otimes B_{j}) = Sp(A_{i}) Sp(B_{j})$; thus
the right hand side of the equality is a morphism from
$\gr(A \otimes B) = A_{0} \otimes B_{0}$ to $A \otimes B$ satisfying
condition (\ref{equation:conditionpolaire}) on poles: it has to be 
$S_{\overline{c}} \hat{F}_{A \otimes B}$. \\
(ii) From the functoriality of the filtration, we know that $F$
only has rectangular blocks relating slopes $\mu_{i} \leq \nu_{j}$,
and that $F_{0}$ is made up of those such that $\mu_{i} = \nu_{j}$.
It is sensible to call the latter ``diagonal blocks''. Then, the
compositum
$(S_{\overline{c}} \hat{F}_{B})^{-1} F \; S_{\overline{c}} \hat{F}_{A}$
is a (meromorphic) morphism from $A_{0}$ to $B_{0}$, with diagonal $F_{0}$
(since $S_{\overline{c}} \hat{F}_{A}$ and $S_{\overline{c}} \hat{F}_{B}$
are in $\G$) and with $0$ under the diagonal. Any block $F_{i,j}$ such that
$\mu_{i} < \nu_{j}$ has all its poles on $[-c;q]$, and with multiplicities
$\leq \nu_{j} - \mu_{i}$; thus, 
$F_{i,j} = \theta_{c}^{\mu_{i} - \nu_{j}} F'_{i,j}$, where $F'_{i,j}$
is holomorphic on $\C^{*}$ and satisfies:
$\sq F'_{i,j} c^{\mu_{j}} A_{j} = c^{\mu_{i}} A_{i} F'_{i,j}$
The same computation (with the Laurent series) as in section
\ref{subsection:Stokesoperators} shows that, under condition
(\ref{equation:nonresonancemixtefonctorielle}), this implies
$F'_{i,j} = 0$. Therefore 
$(S_{\overline{c}} \hat{F}_{B})^{-1} F \; S_{\overline{c}} \hat{F}_{A} = F_{0}$
and (\ref{equation:sommationfonctorielle}) holds.
\hfill $\Box$ \\

In terms of the fiber functors introduced after theorem 
\ref{theorem:filtration}, the meaning of the above lemma is that, 
under proper restrictions to ensure that $S_{\overline{c}} \hat{F}_{A}$ 
is well defined at $a$,and that the nonresonancy conditions 
(\ref{equation:nonresonancemixtetensorielle}) and
(\ref{equation:nonresonancemixtefonctorielle}) hold for any pair of objects, 
$A \leadsto S_{\overline{c}} \hat{F}_{A}(a)$ 
is an $\otimes$-isomorphism 
from $\hat{\omega}^{(0)}_{a}$ to $\omega^{(0)}_{a}$. 
For any pure $A_{0}$, taking up the previous notations, we therefore define,
first its ``weighted spectrum'' and singular locus:
$$
\overline{WSp}(A_{0}) = \text{~the subgroup of $\Eq \times \Z$ generated by~}
\bigcup_{i} \bigl(\overline{Sp(A_{i})} \times \{\mu_{i}\}\bigr),
\quad
\tilde{\Sigma}(A_{0}) = \bigcup_{\mu \not= 0}
\{\overline{c} \in \Eq \;/\; 
(\mu \overline{c},\mu) \in \overline{WSp}(A_{0})\}.
$$
\begin{prop}
Let $<A>$ be the tannakian subcategory of $\EE_{1}^{(0)}$ generated by $A$.
Fix $\overline{c} \not\in \tilde{\Sigma}(A_{0})$ and $a \not\in [-c;q]$.
Then $B \leadsto S_{\overline{c}} \hat{F}_{B}(a)$ is an $\otimes$-isomorphism 
from $\hat{\omega}^{(0)}_{a}$ to $\omega^{(0)}_{a}$, both being restricted 
to $<A>$.
\end{prop}
\Pr
Apply the lemma and the formulas giving the slopes of linear constructions 
in \cite{JSFIL}.
\hfill $\Box$ \\

Now we recall, from \cite{JSSTO} that, for a pure object $A = A_{0}$,
all the $S_{\overline{c}} \hat{F}_{A}$ are equal (they are indeed equal
to the formal Stokes operator $\hat{F}_{A}$ which is actually analytic).

\begin{thm}
With the same restrictions, fix an arbitrary 
$\overline{c_{0}} \not\in \tilde{\Sigma}(A_{0}) \cup \{\overline{-a}\}$.
Then, for all 
$\overline{c} \not\in \tilde{\Sigma}(A_{0}) \cup \{\overline{-a}\}$,
we have, using notation (\ref{equation:defoperateurStokes}):
$$
S_{\overline{c_{0}},\overline{c}} \hat{F}_{A}(a) \in \St(A).
$$
\end{thm}
\Pr
By the above proposition, it is in the Galois group; by the remark above, it is
killed by the functor $\gr$.
\hfill $\Box$ \\

\begin{cor}
We get a family of elements of Lie-like automorphisms:
$\Lsca(A) \underset{def}{=}
\log (S_{\overline{c_{0}},\overline{c}} \hat{F}_{A}(a)) \in \st(A)$.
\end{cor}

Now, although the functoriality and $\otimes$-compatibility were
proved only for 
$\overline{c} \not\in \tilde{\Sigma}(A_{0}) \cup \{\overline{-a}\}$,
the above formula is actually well defined for all
$\overline{c} \not\in \Sigma_{A_{0}} \cup \{\overline{-a}\}$.
Moreover, from the explicit computation in section 
\ref{subsection:Stokesoperators} (multiplications by powers 
of $\theta_{c}$ and resolution of recursive equations by 
inversion of $\Phi_{q^{p} A'_{j},A'_{i}}$), we see that the
mapping $\overline{c} \mapsto \Lsca(A)$ is
meromorphic on $\Eq$, with poles on $\Sigma_{A_{0}}$. 
Moreover, it takes values in the vector space $\st(A)$
for all $\overline{c}$ except for a denumerable subset:
therefore, it takes all its values in $\st(A)$. Last,
taking residues at a pole is an integration process and
gives values in the same vector space.

\begin{thm}
Define the $q$-alien derivations by the formula:
$$
\Derc(A) = Res_{\overline{d} = \overline{c}} LS_{\overline{d},a}(A).
$$
Then, $\Derc(A) \in \st(A)$. (In order to alleviate 
the notation, we do not mention the arbitrary basepoint $a \in \C^{*}$.)
\end{thm}

Of course, for $\overline{c} \not\in \Sigma_{A_{0}}$, we have 
$\Derc(A) = 0$. 
According to the graduation of $\st$ described at the end of section
\ref{subsection:filtration}, each alien derivation admits a canonical
decomposition:
\begin{equation}
\label{equation:decompositionalienderivation}
\Der = \bigoplus_{\delta \geq 1} \Derdc,
\end{equation}
where $\Derdc(A) \in \st_{\delta}(A)$ has only non null blocks for
$\mu_{j} - \mu_{i} = \delta$.

\begin{thm}
The alien derivations are Lie-like $\otimes$-endomorphisms of 
$\hat{\omega}^{(0)}_{a}$ over $\EE_{1}^{(0)}$.
\end{thm}
\Pr
This means first that they are functorial; for all morphisms 
$F : A \rightarrow B$, one has:
$$
\Derc(B) \circ \hat{\omega}^{(0)}_{a}(F) =
\hat{\omega}^{(0)}_{a}(F) \circ \Derc(A).
$$
Note that $\hat{\omega}^{(0)}_{a}(F) = F_{0}(a)$.
First assume the previous restrictions on $\overline{c}$. 
Then, from the lemma \ref{lemma:tensfonctrestreint}, we get that
$S_{\overline{c_{0}},\overline{c}} \hat{F}_{B}(a) \circ F_{0}(a) =
F_{0}(a) \circ S_{\overline{c_{0}},\overline{c}} \hat{F}_{A}(a)$.
Now, the logarithm of a unipotent matrix $P$ being a polynomial
of $P$, we have
$R \circ Q = Q \circ P \Rightarrow \log R \circ Q = Q \circ \log P$,
so that we have
$\Lsca(B) \circ F_{0}(a) = 
F_{0}(a) \circ \Lsca(A)$ and we take the residues on both sides.
Now that the equality is established outside a denumerable set of values of
$\overline{c}$, we can extend it to all values by holomorphy. \\
The assertion means, second, Lie-like $\otimes$-compatibility:
$$
\Derc(A \otimes B) = 
1 \otimes \Derc(B) + \Derc(A) \otimes 1,
$$
where the left and right $1$ respectively denote the identities of
$\hat{\omega}^{(0)}_{a}(A)$ and $\hat{\omega}^{(0)}_{a}(B)$. This
equality makes sense because $\hat{\omega}^{(0)}_{a}$ is itself
$\otimes$-compatible. From the lemma \ref{lemma:tensfonctrestreint}
we get first that
$\Lsca(A \otimes B) = 
\Lsca(A) \otimes \Lsca(A)$. Then, we note
that, for any two unipotent matrices $P$ and $Q$, the commuting product
$P \otimes Q = (P \otimes 1) (1 \otimes Q) = (1 \otimes Q) (P \otimes 1)$
entails $\log (P \otimes Q) = (\log P) \otimes 1 + 1 \otimes (\log Q)$.
The proof is then finished as above.
\hfill $\Box$ 


\subsection{Alien derivations and $q$-Borel transform}

Let $\delta \in \N^{*}$ and assume that the matrix $A$ of 
(\ref{equation:forme-canonique}) has only null blocks $U_{i,j}$ 
for $\mu_{j} - \mu_{i} < \delta$. Then, in the computation
(\ref{equation:calculrecursifhat(F)}), we find the following
equation for $\mu_{j} - \mu_{i} < \delta$:
$\sq \hat{F}_{i,j} z^{\mu_{j}} A_{j} - z^{\mu_{i}} A_{i} \hat{F}_{i,j} = 0$.
Likewise, the upper diagonal blocks of any $S_{\overline{c}} \hat{F}_{A}$
satisfy exactly the same equations. These have no non trivial formal 
solution, neither non trivial meromorphic wih less than $(\mu_{j} - \mu_{i})$ 
poles modulo $q^{\Z}$ (this follows from \ref{subsubsection:caseofonelevel}).
Hence, as well $\hat{F}_{A}$ as all the summations 
$S_{\overline{c}} \hat{F}_{A}$ have null blocks $F_{i,j}$ 
for $0 < \mu_{j} - \mu_{i} < \delta$. \\
On level $\mu_{j} - \mu_{i} = \delta$, the equations to be solved are:
$$
\sq F_{i,j} z^{\mu_{j}} A_{j} - z^{\mu_{i}} A_{i} F_{i,j} = U_{i,j},
$$
which is of the same type as those of \ref{subsubsection:caseofonelevel}.
The properties of this \emph{first non trivial level} of $\hat{F}_{A}$ and
$S_{\overline{c}} \hat{F}_{A}$ will play a crucial role in \cite{RS2}.
Indeed, the logarithm $\log F$ has, as first non trivial level the same
level $\delta$, and the corresponding diagonal is equal to that of $F$.
Therefore, after taking residues, on gets straightaway the $\Derdc$.
To study it in some detail, we therefore take again the light notations
of \ref{subsubsection:caseofonelevel}.

\paragraph{Solving (\ref{equation:onelevelequation1}) 
with $q$-Borel transforms.}
We consider $\delta \in \N^{*}$ as fixed, to alleviate notations.
Let the Laurent series expansion:
$$
\theta^{\delta} = \sum_{n \in \Z} t_{n} z^{n}.
$$
Then, from the functional equation 
$\sq \theta^{\delta} = z^{\delta} \theta^{\delta}$,
we draw the recurrence relations:
$$
\forall n \in \Z \;,\; t_{n-\delta} = q^{n} t_{n}.
$$
>From this, we get the useful estimation:
$$
t_{n} \approx |q|^{- n^{2}/2 \delta}.
$$
The notation $u_{n} \approx v_{n}$ for positive sequences here means 
``same order of magnitude up to a polynomial factor'', more precisely:
$$
u_{n} \approx v_{n} \Longleftrightarrow \exists \; R > 0 \;:\; 
u_{n} = O(R^{n} v_{n}) \text{~and~} v_{n} = O(R^{n} u_{n}).
$$
Now, for any Laurent series 
$F(z) = \sum F_{n} z^{n} \in E \otimes \C[[z,z^{-1}]]$ with coefficients 
$F_{n}$ in some finite dimensional $\C$-vector space $E$, we define its
$q$-Borel transform at level $\delta$ by the formula:
$$
\Bd F(\xi) = \sum t_{-n} F_{n} \xi^{n} \in E \otimes \C[[\xi,\xi^{-1}]].
$$
This transformation strongly increases the convergence properties; for 
instance, if $F \in E \otimes \Ra$, then $\Bd F \in E \otimes \Ree$, etc.
Since we are interested in analyticity of $\Bd F$, we introduce conditions 
on the order of growth of coefficients, adapted from \cite{RamisGrowth}.
Let $G \in E \otimes \C((\xi)) = E \otimes \C[[\xi]][\xi^{-1}]$.
We say that $G := \sum G_{n} \xi^{n} \in E \otimes \Kaqd$ if 
$\parallel G_{n} \parallel = O(R^{n} q^{- n^{2}/2d})$ for some $R > 0$.
We say that $G \in E \otimes \Kaq(d)$ if 
$\parallel G_{n} \parallel = O(R^{n} q^{- n^{2}/2d})$ for all $R > 0$.
In the case that, moreover, $G$ has no pole at $0$ ($G \in \C[[\xi]]$),
we respectively say that $G \in E \otimes \Raqd$,
resp. $G \in E \otimes \Raq(d)$.
Thus, we have the obvious equivalences:
\begin{eqnarray*}
F \in E \otimes \Ra & \Longleftrightarrow & \Bd F \in E \otimes \Raqd, \\
F \in E \otimes \Ka & \Longleftrightarrow & \Bd F \in E \otimes \Kaqd, \\
F \in E \otimes \Ree & \Longleftrightarrow & \Bd F \in E \otimes \Raq(d), \\
F \in E \otimes \Ke & \Longleftrightarrow & \Bd F \in E \otimes \Kaq(d).
\end{eqnarray*}

With the notations of equation (\ref{equation:onelevelequation1}),
write $G = \Bd F = \sum G_{n} \xi^{n}$ and $V = \Bd U = \sum G_{n} \xi^{n}$
(so that $G_{n} = t_{-n} F_{n}$ and $V_{n} = t_{-n} U_{n}$).
Then, multiplying relation (\ref{equation:onelevelequation2}) 
by $t_{-n}$ and noting that $q^{n-\delta} t_{-n} = t_{-(n-\delta)}$, we get:
$$
\forall n \in \Z \;,\; G_{n-\delta} - \Lambda(G_{n}) = V_{n}.
$$
Multiplying by $\xi^{n}$ and summing for $n \in \Z$ yields:
$$
(\xi^{\delta} \Id - \Lambda) \Bd F(\xi) = \Bd U(\xi).
$$
Since $U$ has a positive radius of convergence, $\Bd U(\xi)$ is an entire
function. For $F$ to be a \emph{convergent} solution, it is necessary that
$\Bd F$ be an \emph{entire} function. We shall now appeal to linear algebra.
We first write $A = A_{s} A_{u}$ and $B = B_{s} B_{u}$ the multiplicative
Dunford decompositions. Then $A_{u}^{1/\delta}$ and $B_{u}^{1/\delta}$ are well
defined. In order to define $A_{s}^{1/\delta}$ and $B_{s}^{1/\delta}$, 
it is enough
to choose a mapping $x \mapsto x^{1/\delta}$ on $\C^{*}$ and to apply it
to the eigenvalues. We then put 
$A^{1/\delta} = A_{s}^{1/\delta} A_{u}^{1/\delta}$,
$B^{1/\delta} = B_{s}^{1/\delta} B_{u}^{1/\delta}$ and get a linear map 
$L: F \mapsto A^{1/\delta} F (B^{1/\delta})^{-1}$, which is 
a $\delta^{th}$ root
of $\Lambda$. Call $\mu_{\delta}$ the set of $\delta^{th}$ roots of $1$ 
in $\C$.

\begin{lemma}
Let $E$ be a finite dimensional $\C$-vector space, $A$ an endomorphism
of $E$ and $R$ be any of the following algebras of functions: 
$\Ree$; $\Ree[\xi^{-1}]$; $\Raqd$; $\Kaqd$; $\Raq(d)$; $\Kaq(d)$.
Then the linear operator $(\xi^{\delta} - A^{\delta})$ maps injectively
$E \otimes R$ into itself, its image has a finite codimension $\delta \dim E$
and there is an explicit projection formula on the supplementary
space $E \oplus \cdots \oplus E \xi^{\delta-1}$ of the image:
$$
V \mapsto \sum_{j \in \mu_{\delta}} d (j A)^{\delta-1} P_{j}(A,\xi) V(jA),
$$
where $P_{j}(A,\xi)$ and $V(jA)$ respectively are the following linear 
operator and vector:
$$
P_{j}(A,\xi) = \sum_{i=0}^{\delta-1} (j A)^{i} \, \xi^{\delta-1-i},
\quad
V(jA) = \sum (jA)^{n} V_{n} \in E, \text{~(where~} V = \sum V_{n} \xi^{n}
\text{~is entire)}.
$$
\end{lemma}
\Pr
The algebraic part of the proof rests on the following computation:
$$
1 = \sum_{j \in \mu_{\delta}} \delta (j a)^{\delta-1} P_{j}(a,X), 
\quad \text{~where~}
P_{j}(a,X) = \dfrac{X^{\delta} - a^{\delta}}{X-a} =
\sum_{i=0}^{\delta-1} (j a)^{i} \, X^{\delta-1-i}.
$$
>From this, we draw:
\begin{eqnarray*}
V(\xi) & = & 
\sum_{j \in \mu_{\delta}} \delta (j A)^{\delta-1} P_{j}(A,\xi) V(\xi) \\
& = & \sum_{j \in \mu_{\delta}} 
\delta (j A)^{\delta-1} P_{j}(A,\xi) \bigl(V(\xi) - V(jA)\bigr)
+ \sum_{j \in \mu_{\delta}} \delta (j A)^{\delta-1} P_{j}(A,\xi) V(j A);
\end{eqnarray*}
then we note that, since $P_{j}(A,\xi) (\xi - jA) = \xi^{\delta} - A^{\delta}$,
the first term of the last right hand side is in the image of 
the linear operator $(\xi^{\delta} - A^{\delta})$. The second term plainly
belongs to the supplementary 
space $E \oplus \cdots \oplus E \xi^{\delta-1}$. \\
Then, there are growth conditions on the coefficients to be checked.
In the case of $\Ree$; $\Ree[\xi^{-1}]$, they are standard. In the case
of $\Raqd$, $\Kaqd$, $\Raq(d)$ and $\Kaq(d)$, they follow from the
estimations given in the proof of lemma 2.9 of \cite{JSSTO}.
\hfill $\Box$

\begin{thm}
With the notations of section \ref{subsubsection:caseofonelevel},
equation (\ref{equation:onelevelequation1}) has a convergent solution
if, and only if, $\Bd U(j L) = 0$ for all $j \in \mu_{\delta}$.
More precisely, the family 
$\bigl(\Bd U(j L) \bigr)_{j \in \mu_{\delta}} \in 
\Ma_{r,s}(\C)^{\mu_{\delta}} \simeq \Ma_{r,s}(\C)^{\delta}$ 
is a complete set of invariants for analytic classification within the 
formal class $M_{0}$.
\end{thm}

\paragraph{Solving (\ref{equation:onelevelequation1}) 
with $\theta$ functions and residue invariants.}

>From the computations in the one level case of 
\ref{subsubsection:caseofonelevel},
we see that the only solution of (\ref{equation:onelevelequation1}) 
such that condition (\ref{equation:conditionpolaire}) holds is given 
by the explicit formula:
$$
F_{\overline{c}}(z) = \dfrac{1}{\theta_{c}^{\delta}} 
\sum_{n \in \Z} (c^{\delta} q^{n} \Id - \Lambda)^{-1} \, V_{n} \, z^{n},
\text{~where~} V = \theta_{c}^{\delta} U.
$$
A short computation shows that $V = \sum t_{p} c^{-p} U_{n} z^{n+p}$,
so that, $V_{0} = \sum t_{-n} c^{n} U_{n} = \Bd U(c)$. \\
In order to compute explicitly the alien derivation in the one level 
case, it is convenient to normalize the setting, by requiring that 
all eigenvalues of $A$ and $B$ lie in the fundamental annulus 
$1 \leq |z| < |q|$ (up to shearing transformation, this is always 
possible). We may further decompose the pure blocks $z^{\mu} A$
and $z^{\nu} B$ into their corresponding characteristic subspaces.
In other words, we may (and shall) assume here that $A$ and $B$
are block diagonal, each block $A_{\alpha}$ (resp. $B_{\beta}$)
have the unique eigenvalue $\alpha$ (resp. $\beta$), this lying
in the fundamental annulus. We write $\Lambda_{\alpha,\beta}$,
$L_{\alpha,\beta}$, $U_{\alpha,\beta}$, and compute the corresponding 
component $\Derdxab(M)$ of $\Derdx(M)$.
Let $\xi \in \C^{*}$ be a prohibited (polar) value of $c$. This means that 
one of the matrices $(\xi^{\delta} q^{n} \Id - \Lambda_{\alpha,\beta})$ 
is singular, so
that $\xi^{\delta} q^{n} = \alpha/\beta$. From the normalisation condition, 
we see that this can occur only for one value of $n$. Since residues are 
actually defined on $\Eq$, one can choose $\xi$ such that the bad value 
of $n$ is $n = 0$. Then, we are to compute:
$$
\Derdxab(M) = Res_{c = \xi} \, \dfrac{1}{\theta_{c}^{\delta}(a)} 
(c^{\delta} \Id - \Lambda_{\alpha,\beta})^{-1} \Bd U_{\alpha,\beta}(c).
$$
Note that the arbitrary basepoint $a \in \C^{*}$ (which provides us with
the fiber functor $\hat{\omega}^{(0)}_{a}$) appears only in the theta factor. 
As in the previous section, we introduce $L_{\alpha,\beta}$ such that 
$L_{\alpha,\beta}^{\delta} = \Lambda_{\alpha,\beta}$
and get, from the same formulas as before:
$$
\Derdxab(M) = Res_{c = \xi} \, \dfrac{1}{\theta_{c}^{\delta}(a)} 
\sum_{j \in \mu_{\delta}} \delta (j _{\alpha,\beta}L)^{\delta-1} 
(c - j L_{\alpha,\beta})^{-1} \Bd U_{\alpha,\beta}(c).
$$
Now, $\xi$ is an eigenvalue of one and only one of the $j L_{\alpha,\beta}$,
call it $L_{\xi\alpha,\beta}$. From classical 
``holomorphic functional calculus'' 
(see \emph{e.g.} \cite{Rudin}), we get:
$$
\Derdxab(M) = 
\theta^{-\delta}(L_{\xi\alpha,\beta}^{-1}) \, 
\delta L_{\xi\alpha,\beta}^{\delta-1} \,
\Bd U_{\alpha,\beta}(L_{\xi\alpha,\beta}).
$$
Recall that, as in \emph{loc. cit.}, the theta factor is the application 
of a holomorphic function to a linear operator.

\begin{thm}
Call $\Phi_{a}$ the automorphism of $\Ma_{r,s}(\C)^{\mu_{\delta}}$,
which, on the $(\alpha,\beta)$ component, is left multiplication by
$\theta^{-\delta}(L_{\xi\alpha,\beta}^{-1})$. Then $\Phi_{a}$ sends
the $q$-Borel invariant $\bigl(\Bd U(j L) \bigr)_{j \in \mu_{\delta}}$
to the $\Der$ invariant: $\bigoplus\limits_{\xi} \Derdx(M)$.
\end{thm}

Here is an example similar to that at the end of section 
\ref{subsection:generalfacts}. We take 
$M = \begin{pmatrix} \alpha & u \\ 0 & \beta z \end{pmatrix}$, where 
$\alpha,\beta \in \C^{*}$ and $u \in \Ra$. The slopes are $\mu = 0$ and 
$\nu = 1$ and the only level is $\delta = 1$. The associated
non homogeneous equation is $\beta z \sq f - \alpha f = u$, which, in
the Borel plane, becomes 
$(\beta \xi - \alpha) \mathcal{B}_{q}^{(1)} f = \mathcal{B}_{q}^{(1)} u$,
and the obstruction to finding an analytical solution is the
complex number $\mathcal{B}_{q}^{(1)} u(\alpha/\beta)$. This is also 
the invariant associated to the analytical class of $M$ within 
its formal class. \\
On the side of resolution with $\theta$ and residues, we first get: 
$f_{\overline{c}}(z) = \dfrac{1}{\theta_{c}} 
\sum_{n \in \Z} (c q^{n} - \alpha \beta^{-1})^{-1} \, v_{n} \, z^{n}$,
where $v = \theta_{c} u$, then the only non trivial alien derivation, 
given for $\xi = \alpha/\beta$:
$$
\Derdx = Res_{c = \xi} f_{\overline{c}}(a) =
Res_{c = \xi} \dfrac{1}{\theta_{c}(a)} (c - \xi)^{-1} 
\mathcal{B}_{q}^{(1)} u(c) = 
\dfrac{1}{\theta(a/\xi)} \mathcal{B}_{q}^{(1)} u(\xi).
$$



\section{Conclusion}

The construction of the alien derivations in the differential 
and $q$-difference case are apparently quite different. In fact,
it is possible to reformulate things in the differential case 
to exhibit some analogy; one can mimick the constructions of 
the $q$-difference case: in place of a meromorphic function of 
a $q$-direction of summation in $E_q$,  one gets in the differential 
case a locally constant function of a direction of summation in $S^1$ 
minus a finite singular set, the poles being replaced by ``jump points''. 
The jumps are evaluated by a non-abelian boundary value: one gets 
the Stokes operators. The alien derivation are the logarithms of 
these operators. There is a slight difference with the $q$-difference case: 
in this last case, we took the logarithm \emph{before} evaluating 
the singularity by a residue. We remark that to consider locally 
constant functions on $S^1$ with a finite set of jumps as 
the differential analog of meromorphic functions in the $q$-difference 
case is in perfect accordance with the study of the confluence process 
by the second author in \cite{JSAIF}. \\

Our constructions suggest some interesting problems. \\
1. If we consider the computation of the $q$-alien derivations 
$\Der_{\xi} $ in simple cases, there appears theta factors and factors 
coming from a $q$-Borel transform. In the simplest cases we can 
define (pointed) alien derivations $\Der_{\xi} $ as operators acting 
on some $q$-holonomic power series, and, eliminating the theta factors, 
we can observe that this modified $q$-derivative of a power series 
is itself a power series: we get a new operator, an \emph{unpointed} alien 
derivation $\Delta_{\xi}$. This suggests the possibility to copy the 
Ecalle's definition of alien derivations (cf. \cite{CNPham}) in the 
$q$-difference case: a resurgence lattice in the Borel plane is replaced
by the set of singularities in the different $q$-Borel planes 
corresponding to the different $q$-levels $\delta \in \N^*$ (the $q$-direction 
of summation being fixed), the Ecalle's analytic continuation paths by 
summation paths ``between the levels" and the boundary values by residues. 
In this program it is important to remark that the ``algebraic" definition of summability 
used in this paper is equivalent to ``analytic" definitions in Borel-Laplace style \cite{MZ, RZ, Zha02} \\
2. The global classification: we must put together the work 
of \cite{JSGAL} and the results of the present article (at zero 
and infinity). It is not difficult to guess what will happen and 
to describe a ``fundamental group'' which is Zariski dense in the 
tannakian Galois group, but some great problems remain: even in 
the regular singular case, we know neither the structure 
of the global tannakian Galois group (except in the \emph{abelian} 
case), nor if there exists a reasonable ``localisation 
theory" for the singularities on $\C^*$ (between $0$ and $\infty$). \\
3. The confluence problem. Some simple examples suggest an extension of 
the results of the second author (cf. \cite{JSAIF}) to the irregular case: 
confluence of $q$-Stokes phenomena at $0$ to Stokes phenomena. It is 
natural to study what will happen with the alien derivations (there 
is no hope with the pointed $q$-alien derivations, due to the bad 
properties of $theta$ functions in the confluence processes, but 
it could work nicely with the unpointed $q$-alien derivations). 

\vskip 10pt

\noindent{\large\bf Acknowledgements}

\vskip 5pt

The work of the first author has been partially supported by the NEST EU Math. Project GIFT, FP6-005006-2.



\end{document}